\documentclass[aop]{imsart}

\RequirePackage[OT1]{fontenc}
\RequirePackage{amsthm,amsmath}
\RequirePackage[numbers]{natbib}
\RequirePackage[colorlinks,citecolor=blue,urlcolor=blue]{hyperref}

\startlocaldefs

\usepackage{amssymb}
\usepackage{mathrsfs} 
\usepackage{amsmath}

\usepackage{amsthm}
\usepackage{amsfonts}
\usepackage{graphicx}

\def\thesection{\arabic{section}}
\renewcommand{\theequation}{\thesection.\arabic{equation}}

\newtheorem{theorem}{Theorem}[section]
\newtheorem{lemma}[theorem]{Lemma}
\newtheorem{proposition}[theorem]{Proposition}
\newtheorem{corollary}[theorem]{Corollary}

\newtheorem{fact}[theorem]{Fact}

\newtheorem{defin}[theorem]{Definition}

\theoremstyle{definition}   
\newtheorem{remark}[theorem]{Remark}

\newcommand{\eqnsection}{
\renewcommand{\theequation}{\thesection.\arabic{equation}}
    \makeatletter
    \csname  @addtoreset\endcsname{equation}{section}
    \makeatother}
\eqnsection

\def\r{{\mathbb R}}

\def\P{{\bf P}}
\def\E{{\bf E}}

\def\z{{\mathbb Z}}

\def\ee{\mathrm{e}}
\def\d{\, \mathrm{d}}
\def\T{{\mathbb T}}

\endlocaldefs

\begin{document}

\begin{frontmatter}
\title{The Derrida--Retaux conjecture on recursive models}
\runtitle{The Derrida--Retaux conjecture}

\begin{aug}
\author{\fnms{Xinxing} \snm{Chen}\thanksref{t1}\ead[label=e1]{chenxinx@sjtu.edu.cn}},
\author{\fnms{Victor} \snm{Dagard}\thanksref{t2}\ead[label=e2]{victor.dagard@lps.ens.fr}},
\author{\fnms{Bernard} \snm{Derrida}\thanksref{t2}\ead[label=e3]{derrida@lps.ens.fr}},\\
\author{\fnms{Yueyun} \snm{Hu}\thanksref{t2, t3}\ead[label=e4]{yueyun@math.univ-paris13.fr}},
\author{\fnms{Mikhail} \snm{Lifshits}\thanksref{t4}\ead[label=e5]{mikhail@lifshits.org}},
\and
\author{\fnms{Zhan} \snm{Shi}\thanksref{t2}
\ead[label=e6]{zhan.shi@upmc.fr}
}

\thankstext{t1}{Partially supported by NSFC grants 11771286 and 11531001.}
\thankstext{t2}{Partially supported by ANR project MALIN 16-CE93-0003.}
\thankstext{t3}{Partially supported by ANR project SWIWS 17-CE40-0032-02.}
\thankstext{t4}{Partially supported by  RFBR-DFG grant 20-51-12004.}
\runauthor{X. Chen et al.}

\affiliation{Shanghai Jiaotong University, \'Ecole Normale Sup\'erieure, Coll\`ege de France,   
Universit\'e Sorbonne Paris Nord,  St.~Petersburg State University, and Sorbonne Universit\'e Paris VI}

\address{School of Mathematical Sciences\\
Shanghai Jiaotong University\\
200240 Shanghai, China\\
\printead{e1}\\}

\address{Laboratoire de Physique\\
\'Ecole Normale Sup\'erieure\\
  PSL, CNRS, F-75005 Paris, France\\
\printead{e2}\\}

\address{Coll\`ege de France, PSL\\
11 place Marcelin Berthelot \\
F-75231 Paris Cedex 05, France\\
\printead{e3}\\}

\address{LAGA, USPN\\
99 av.~J-B Cl\'ement \\
F-93430 Villetaneuse, France\\
\printead{e4}\\}

\address{St.~Petersburg State University\\
Department of Mathematics and Computer Sciences\\
199034, St.~Petersburg\\
Universitetskaya emb., 7-9, Russian Federation\\
\printead{e5}\\}

\address{LPSM\\
Sorbonne Universit\'e Paris VI \\
 4 place Jussieu, F-75252 Paris Cedex 05\\
 France\\
\printead{e6}\\}
\end{aug}

\begin{abstract} We are interested in the nearly supercritical regime in a family of max-type recursive models studied by Collet, Eckman, Glaser and Martin~\cite{collet-eckmann-glaser-martin} and by Derrida and Retaux~\cite{derrida-retaux}, and prove that under a suitable integrability assumption on the initial distribution, the free energy vanishes at the transition with an essential singularity with exponent $\tfrac12$. This gives a weaker answer to a conjecture of Derrida and Retaux~\cite{derrida-retaux}. Other behaviours are obtained when the integrability condition is not satisfied.
\end{abstract}



\end{frontmatter}

\section{Introduction}
\label{s:intro}

\subsection{The model and main results}
 
Let $m\ge 2$ be an integer. Let $X_0 \ge 0$ be a random variable taking values in $\z_+ := \{ 0, \, 1, \, 2, \ldots\}$; to avoid triviality, it is assumed, throughout the paper, that
$\P(X_0\ge 2)>0$. Consider the following recurrence relation: for all $n\ge 0$,
\begin{equation}
    X_{n+1}
    =
    (X_{n,1} + \cdots + X_{n,m} -1)^+ ,
    \label{iteration}
\end{equation}

\noindent where $X_{n,i}$, $i\ge 1$, are independent copies of $X_n$. Notation: $x^+ := \max\{ x, \, 0\}$ for all $x\in \r$.

{F}rom \eqref{iteration}, we get $m \, \E(X_n) -1 \le \E(X_{n+1}) \le m \, \E(X_n)$, which enables us to define the {\it free energy}
\begin{equation}
    F_\infty
    :=
    \lim_{n\to \infty} \downarrow \, \frac{\E(X_n)}{m^n}
    =
    \lim_{n\to \infty} \uparrow \, \frac{\E(X_n) - \frac{1}{m-1}}{m^n}
    \ge
    0 \, .
    \label{F}
\end{equation}

We now recall a conjecture of Derrida and Retaux~\cite{derrida-retaux}. For any random variable $X$, we write $P_X$ for its law. Assume
$$
P_{X_0}
=
(1-p) \, \delta_0 + p\, P_{X_0^*},
$$

\noindent where $\delta_0$ denotes the Dirac measure at $0$, $X_0^*$ a (strictly) positive integer-valued random variable satisfying $\P(X_0^*\ge 2)>0$, and $p\in [0, \, 1]$ a parameter. Since $p\mapsto F_\infty =: F_\infty(p)$ is non-decreasing, there exists $p_c = p_c(X_0^*) \in [0, \, 1]$ such that $F_\infty (p)>0$ for $p>p_c$ and that $F_\infty (p) =0$ for $p<p_c$.\footnote{We are going to see that $p_c<1$.} The Derrida--Retaux conjecture says that if $p_c>0$ (and possibly under some additional integrability conditions on $X_0$), then
\begin{equation}
     F_\infty(p)
     =
     \exp \Big( - \frac{C+o(1)}{(p-p_c)^{1/2}} \Big),
     \qquad
     p\downarrow p_c\, ,
     \label{conj:bernard}
\end{equation}

\noindent for some constant $C\in (0, \, \infty)$. When $p_c=0$, it is possible to have other exponents than $\frac12$ in \eqref{conj:bernard}, see \cite{yz_bnyz}. In \cite{bmxyz_questions}, we have presented several open questions concerning the critical regime $p=p_c$ when $p_c>0$.

The model with recursion defined in \eqref{iteration}  is known to have a phase transition (Collet et al.~\cite{collet-eckmann-glaser-martin}), recalled in Theorem A below. It is expected to have many universal properties at or near criticality, though few of these predicted properties have been rigorously proved so far. The model was introduced by Derrida and Retaux~\cite{derrida-retaux} as a simplified hierarchical renormalization model to understand the depinning transition of a line in presence of strong disorder \cite{derrida-hakim-vannimenus}. The exponent $\frac12$ in the Derrida--Retaux conjecture \eqref{conj:bernard} was already predicted in \cite{tang-charte}, while another exponent $1$ was predicted in \cite{monthus}. For the mathematical literature of pinning models, see \cite{alexander}, \cite{giacomin}, \cite{giacomin_stf}, \cite{giacomin-toninelli}, \cite{giacomin-toninelli-lacoin}. The same exponent 1 was found for a copolymer model in \cite{berger-giacomin-lacoin} with additional precision. The recursion \eqref{iteration} has appeared in another context, as a spin glass toy-model in Collet et al.\ \cite{collet-eckmann-glaser-martin2}--\cite{collet-eckmann-glaser-martin}, and is moreover connected to a parking scheme investigated by Goldschmidt and Przykucki~\cite{goldschmidt-przykucki}; it was studied from the point of view of iterations of random functions (Li and Rogers~\cite{li-rogers}, Jordan~\cite{jordan}), and also figured as a special case in the family of max-type recursive models analyzed in the seminal paper of Aldous and Bandyopadhyay~\cite{aldous-bandyopadhyay}. See Hu and Shi~\cite{yz_bnyz} for an extension to the case when $m$ is random, and Hu, Mallein and Pain~\cite{HMP} for an exactly solvable version in continuous time.

The aim of this paper is to study the Derrida--Retaux conjecture. Let us first recall the following characterisation of the critical regime.

\bigskip

\noindent {\bf Theorem A (Collet et al.~\cite{collet-eckmann-glaser-martin}).} {\it We have
\begin{equation}
    p_c(X_0^*)
    =
    \frac{1}{1+ \E\{ [(m-1)X_0^*-1] m^{X_0^*}\} }
    \in [0, \, 1)\, .
    \label{p_c}
\end{equation}

\noindent In words, assuming $\E(X_0^* \, m^{X_0^*}) <\infty$, then $p=p_c$ means $\E(m^{X_0}) = (m-1) \E(X_0 \, m^{X_0})$, and more precisely, $F_\infty>0$ if $p>p_c$, while $F_\infty =0$ if $p\le p_c$.}

\bigskip

It is natural to say that the system is subcritical if $p<p_c$, critical if $p=p_c$, and supercritical if $p>p_c$. Note from Theorem A that the assumption $p_c>0$ in the Derrida--Retaux conjecture is equivalent to saying that $\E(X_0^* \, m^{X_0^*}) <\infty$.

We give a partial answer to the Derrida--Retaux conjecture, by showing that under suitable general assumptions on the initial distribution, $\frac12$ is the correct exponent, in the exponential scale, for the free energy.

\medskip

\begin{theorem}
\label{t:main}

 Assume $\E[(X_0^*)^3 \, m^{X_0^*}] <\infty$. Then
 $$
   F_\infty (p)
   =
    \exp \Big( - \frac{1}{(p-p_c)^{\frac12 + o(1)}} \Big) ,
   \qquad
   p \downarrow p_c\, .
 $$

\end{theorem}

\medskip

It is possible to  obtain some information  about $o(1)$ in Theorem \ref{t:main}; see \eqref{ub} and \eqref{lb}. A similar remark applies to Theorem \ref{t:2<alpha<4} below.

It turns out that our argument in the proof of the lower bound (for the free energy) in Theorem \ref{t:main} is quite robust. With some additional minor effort, it can be adapted to deal with systems that do not satisfy the condition $\E[(X_0^*)^3 \, m^{X_0^*}] <\infty$. Although this integrability condition might look exotic, it is optimal for the validity of the Derrida--Retaux conjecture. In the next theorem, we assume\footnote{Notation: by $a_k \sim b_k$, $k\to \infty$, we mean $\lim_{k\to \infty} \frac{a_k}{b_k} =1$.} $\P(X_0^* = k) \sim c_0\, m^{-k} k^{-\alpha}$, $k\to \infty$, for some constant $0<c_0<\infty$ and some parameter $2<\alpha \le 4$. The inequality $\alpha>2$ ensures $p_c>0$, which is the basic condition in the Derrida--Retaux conjecture, whereas the inequality $\alpha\le 4$ implies $\E[(X_0^*)^3 \, m^{X_0^*}] =\infty$. It turns out that in this case, the behaviour of the free energy differs from the prediction in the Derrida--Retaux conjecture.

\medskip

\begin{theorem}
\label{t:2<alpha<4}

 Assume $\P(X_0^* = k) \sim c_0\, m^{-k} k^{-\alpha}$, $k\to \infty$, for some $0<c_0<\infty$ and $2<\alpha \le 4$. Then
 $$
 F_\infty(p)
 =
 \exp \Big( - \frac{1}{(p-p_c)^{\nu + o(1)}} \Big) ,
 \qquad
 p \downarrow p_c\, ,
 $$
 where $\nu = \nu(\alpha) := \frac{1}{\alpha-2}$.
\end{theorem}

\medskip

Let us keep considering the situation $\P(X_0^* = k) \sim c_0\, m^{-k} k^{-\alpha}$, $k\to \infty$, for some $0<c_0<\infty$. The case $2<\alpha\le 4$ was considered in Theorem \ref{t:2<alpha<4}. When $-\infty<\alpha \le 2$, we have $\E(X_0^* \, m^{X_0^*}) = \infty$, which violates the basic condition $p_c>0$ in the Derrida--Retaux conjecture; so the conjecture does not apply to this situation. In \cite{yz_bnyz}, it was proved that if $-\infty<\alpha<2$, then $F_\infty (p) = \exp ( - \frac{1}{p^{\nu_1 + o(1)}})$ when $p\downarrow p_c=0$, with $\nu_1 = \nu_1 (\alpha) = \frac{1}{2-\alpha}$. This leaves us with the case $\alpha=2$.

\medskip

\begin{theorem}
\label{t:alpha=2}

 Assume $\P(X_0^* = k) \sim c_0\, m^{-k} k^{-2}$, $k\to \infty$, for some $0<c_0<\infty$. Then
 $$
 F_\infty (p)
 =
 \exp \Big( - \ee^{(C+o(1))/p}\Big) ,
 \qquad
 p \downarrow p_c=0\, ,
 $$
 where $C:= \frac{1}{(m-1)c_0}$.
\end{theorem}


\subsection{Description of the proof}

Having in mind both the supercritical system (in Theorem \ref{t:main}) and the system with initial distribution satisfying $\P(X_0^* = k) \sim c_0\, m^{-k} k^{-\alpha}$ with $2<\alpha\le 4$ or $\alpha=2$ (in Theorems \ref{t:2<alpha<4} and \ref{t:alpha=2}, respectively), we introduce in Section \ref{s:regularity} a notion of {\bf regularity} for systems. It is immediately seen that supercritical systems are regular (Lemma \ref{A_4P(X_0=0)}), and so are appropriately truncated systems in Theorems \ref{t:2<alpha<4} and \ref{t:alpha=2} (Lemma \ref{e:chi_lower_bound} in Section \ref{subs:light_tail_lb}). Most of forthcoming technical results are formulated for regular systems in view of applications in the proof of Theorem \ref{t:main} on the one hand, and of Theorems \ref{t:2<alpha<4} and \ref{t:alpha=2} on the other hand.

To study the free energy $F_\infty$, we make the simple observation that by \eqref{F}, for all $n\ge 0$,
\begin{equation}
    \frac{\E(X_n) - \frac{1}{m-1}}{m^n}
    \le
    F_\infty
    \le
    \frac{\E(X_n)}{m^n} \, ,
    \qquad
    n\ge 0\, .
    \label{F_encadrement}
\end{equation}

\noindent So in order to bound $F_\infty$ from above, we only need to find a sufficiently large $n$ such that $\E(X_n)\le 3$ (say), whereas to bound $F_\infty$ from below, it suffices to find another $n$ not too large, for which $\E(X_n) \ge 2$.

Upper bound: the upper bound in the theorems is proved by studying the moment generating function. In the literature, the moment generating function is a commonly used tool to study the recursive system (\cite{collet-eckmann-glaser-martin}, \cite{derrida-retaux}, \cite{bmxyz_questions}). In Section \ref{s:ub}, we obtain a general upper bound (Proposition \ref{p:ub}) for $\E(X_n)$ for regular systems. Applying Proposition \ref{p:ub} to the supercritical system yields the upper bound in Theorem \ref{t:main}. 

Lower bound: the proof of the lower bound in Theorem \ref{t:main} requires some preparation. In Section \ref{s:preliminaries}, an elementary coupling, called the $XY$-coupling in Theorem \ref{t:coupling}, is presented for the supercritical system $(X_n)$ and a critical system $(Y_n)$, in such a way that $Y_n \le X_n$ for all $n$. We then use a natural hierarchical representation of the systems, and study $N_n^{(0)}$, the number of open paths (the paths on the genealogical tree along which the operation $x\mapsto x^+$ is unnecessary) up to generation $n$ with initial zero value. The most important result in Section \ref{s:preliminaries} is the following inequality: if $\E(X_0-Y_0) \ge \eta \, \P(Y_0=0)$ for some $\eta>0$, then for suitable non-negative $r$, $n$, $k$ and $\ell$,
$$
     \E(X_{n+k+\ell})
     \ge
     \frac{m^{k+\ell} \, \eta}{2}\, \E\Big[ N_n^{(0)}\, {\bf 1}_{\{ N_n^{(0)} \ge r\} } \, {\bf 1}_{\{ Y_n =k\} } \Big] \, ;
     \leqno(\ref{pf_lb_E(Xn)>_eq3})
$$

\noindent see Theorem \ref{t:connexion}. This inequality serves as a bridge connecting, on the one hand, the expected value of the supercritical system $(X_n)$, and on the other hand, the expected number of open paths in the critical system $(Y_n)$. 

Using \eqref{pf_lb_E(Xn)>_eq3} and the upper bound for $\E(X_n)$ established in Proposition \ref{p:ub}, we obtain an upper bound for $\E [ m^{Y_n} \, N_n^{(0)}\, {\bf 1}_{\{ N_n^{(0)} \ge r\} } \, {\bf 1}_{\{ Y_n \le n\} } ]$ for all $n$ and suitable $r = r(n)$; see Corollary \ref{c:E(N)>}. This application of \eqref{pf_lb_E(Xn)>_eq3} is referred to as the first crossing of the bridge, and is relatively effortless. 
 
We intend to cross the bridge for a second time, but in the {\it opposite} direction. To prepare for the second crossing, we prove a recursive formula for $\E[m^{Y_n}(1+Y_n)N_n^{(0)}]$ in Proposition \ref{p:recursion_number_open_paths}:
$$
\E[m^{Y_n}(1+Y_n)N_n^{(0)}]
=
\P(Y_0=0) \, \prod_{k=0}^{n-1} [\E(m^{Y_k})]^{m-1} \, .
$$

\noindent Since the asymptotics of $\prod_{k=0}^{n-1} [\E(m^{Y_k})]^{m-1}$ are known, this formula gives useful upper and lower bounds for $\E[m^{Y_n}(1+Y_n)N_n^{(0)}]$, stated in \eqref{joint_moment_Yn_Nn}. 

We are now ready to establish a good upper bound for $\E(m^{Y_n}N_n^{(0)})$ (see Lemma \ref{sum_I_n}): on the event $\{N_n^{(0)} \ge r, \, Y_n \le n \}$, the expectation was already handled  by  Corollary \ref{c:E(N)>}; on the complementary of this event (which needs to be split into two sub-cases), an application of the Markov inequality does the job thanks to the upper bound for $\E[m^{Y_n}(1+Y_n)N_n^{(0)}]$ in \eqref{joint_moment_Yn_Nn}. [Actually Lemma \ref{sum_I_n} states slightly less: it gives an upper bound only for the Ces\`aro sum of $\E(m^{Y_n}N_n^{(0)})$ --- which nonetheless suffices for our needs.] 

The next step is to write a recursion formula for $\E[(1+Y_n)^3 m^{Y_n} N_n^{(0)}]$ in the same spirit as Proposition \ref{p:recursion_number_open_paths}; together with the upper bound for $\E(m^{Y_n}N_n^{(0)})$ in Lemma \ref{sum_I_n}, the formula gives an upper bound for $\E[(1+Y_n)^3 m^{Y_n} N_n^{(0)}]$: this is Proposition \ref{p:third_derivative_In}. [The power $3$ in $(1+Y_n)$ is important as we are going to see soon.]  

Let us write, for positive integers $r$ and $K$,
\begin{eqnarray*}
 &&\E [ (1+Y_n) m^{Y_n} N_n^{(0)} \, {\bf 1}_{\{ N_n^{(0)} \ge r\} } \, {\bf 1}_{\{ Y_n <K \} } ]
    \\
 &\ge& \E [ (1+Y_n) m^{Y_n} N_n^{(0)}]
    -
    \E [ (1+Y_n) m^{Y_n}\, N_n^{(0)} \, {\bf 1}_{\{ Y_n \ge K \} } ]
    \\
 && \qquad\qquad
    -
    r\, \E [ (1+Y_n) m^{Y_n}]\, .
\end{eqnarray*}

\noindent We can bound $\E [ (1+Y_n) m^{Y_n} N_n^{(0)}]$ from below by means of Proposition \ref{p:recursion_number_open_paths} (or rather its consequence \eqref{joint_moment_Yn_Nn}), and bound $\E [ (1+Y_n) m^{Y_n}\, N_n^{(0)} \, {\bf 1}_{\{ Y_n \ge K \} } ]$ from above by the Markov inequality and Proposition \ref{p:third_derivative_In} (which is why the factor $(1+Y_n)^3$ in the proposition is important, otherwise the bound would not be good enough), whereas $\E [ (1+Y_n) m^{Y_n}]$ is smaller than a constant depending only on $m$. Consequently,   we can choose appropriate values for $r$ and $K$ (both depending on $n$) and obtain a lower bound for $\E [ (1+Y_n) m^{Y_n}\, N_n^{(0)} \, {\bf 1}_{\{ N_n^{(0)} \ge r\} } \, {\bf 1}_{\{ Y_n <K \} } ]$. Since
\begin{eqnarray*}
 &&\max_{k\in [0, \, K)\cap \z} (1+k) m^k \E [ N_n^{(0)} \, {\bf 1}_{\{ N_n^{(0)} \ge r\} } \, {\bf 1}_{\{ Y_n = k \} } ]
    \\
 &\ge& \frac1K \, \E [ (1+Y_n) m^{Y_n}\, N_n^{(0)} \, {\bf 1}_{\{ N_n^{(0)} \ge r\} } \, {\bf 1}_{\{ Y_n <K \} } ] \, ,
\end{eqnarray*}

\noindent there exists an integer $k \in [0, \, K)$ for which we have a lower bound for $\E [ N_n^{(0)} \, {\bf 1}_{\{ N_n^{(0)} \ge r\} } \, {\bf 1}_{\{ Y_n = k \} } ]$. The parameters are chosen such that when using the bridge inequality \eqref{pf_lb_E(Xn)>_eq3} for the second time, we get $\E(X_{n+k+\ell}) \ge 2$ for convenient $n$, $k$ and $\ell$. Together with the first inequality in \eqref{F_encadrement}, this yields the lower bound in Theorem \ref{t:main}. 

Finally in Section \ref{s:alpha}, we prove Theorems \ref{t:2<alpha<4} and \ref{t:alpha=2} by means of a truncation  argument.

The rest of the paper is as follows: 


$\bullet$ Section \ref{s:regularity}: a notion of regularity for systems;

$\bullet$ Section \ref{s:ub}: proof of the upper bound in Theorem \ref{t:main};

$\bullet$ Section \ref{s:preliminaries}: $XY$-coupling, open paths, a lower bound via open paths;

$\bullet$ Section \ref{s:open_paths}: a formula for the number of open paths and other preparatory work;

$\bullet$ Section \ref{s:2nd_crossing}: a general lower bound for free energy;

$\bullet$ Section \ref{s:lb}: proof of the lower bound in Theorem \ref{t:main};

$\bullet$ Section \ref{s:alpha}: proof of Theorems \ref{t:2<alpha<4} and \ref{t:alpha=2}.

\medskip

Notation: we often write $\P_{\! p}$ instead of $\P$ in order to stress dependence on the parameter $p$ (accordingly, the corresponding expectation is denoted by $\E_p$), and $\P_{\! p_c+\varepsilon}$ if $p=p_c+\varepsilon$. When we take $p=p_c+\varepsilon$, it is implicitly assumed that $\varepsilon \in (0, \, 1-p_c)$. Moreover, $\frac{a}{\Lambda} := 0$ for $\Lambda=\infty$ and $a\in \r$.

\section{A notion of regularity for systems}
\label{s:regularity}

Consider a generic $\z_+$-valued system $(X_n, \, n\ge 0)$ defined by \eqref{iteration}, with $\E(X_0 \, m^{X_0}) < \infty$ (which is equivalent to $p_c>0$) and $\P(X_0 \ge 2)>0$.

In several situations we will assume that the system   satisfies a certain additional regularity condition, which is, actually, satisfied if the system is critical or supercritical, or if it is suitably truncated.
\medskip

\begin{defin}

 Let $\zeta$ be a $\z_+$-valued random variable, with $\E(\zeta \, m^\zeta) < \infty$ and $\P(\zeta\ge 2)>0$. Let $\beta\in [0, \, 2]$ and $\chi \in (0, \, 1]$.
 Write
\begin{eqnarray}
    \Lambda (\zeta)
 &:=& \E(\zeta^3\, m^\zeta) \in (0, \, \infty] \, ,
    \label{Lambda}
    \\
    \Xi_k(\zeta)
 &:=& \E \Big[ (\zeta \wedge k)^2 \, [(m-1)\zeta-1] \, m^\zeta \Big] \in (0, \, \infty),
     \qquad
     k\ge 1\, ,
     \label{Theta_0}
\end{eqnarray}

\noindent where $a\wedge b := \min\{ a, \, b\}$.

We say that the random variable $\zeta$ is $\beta$-regular with coefficient $\chi$ if for all integers $k\ge 1$,
 \begin{equation}
     \Xi_k(\zeta)
     \ge
     \chi\, \min\{ \Lambda (\zeta), \, k^\beta\} \, .
     \label{alpha}
 \end{equation}

Furthermore, we say that a system $(X_n, \, n\ge 0)$ is $\beta$-regular with coefficient $\chi$, if $X_0$ is $\beta$-regular with coefficient $\chi$.

\end{defin}

\medskip

\begin{remark}
\label{r:alpha=2;Lambda<infty}

It is immediately seen that if the $\z_+$-valued random variable $\zeta$ is $2$-regular, then $\Lambda (\zeta) <\infty$ (otherwise, \eqref{alpha} would become: $\Xi_k(\zeta) \ge \chi\, k^2$ for all integers $k\ge 1$, which would lead to a contradiction, because $\lim_{k\to \infty} \frac{\Xi_k(\zeta)}{k^2} =0$ by the dominated convergence theorem).\qed

\end{remark}

\begin{remark}

Let us say a few words about our interest in the notion of the system's regularity. Part of our concern is to obtain an upper bound for the moment generating function (the forthcoming Proposition \ref{p:ub}, for supercritical systems), which plays a crucial role in {\it both} upper and lower bounds in Theorems \ref{t:main}, \ref{t:2<alpha<4} and \ref{t:alpha=2}. In this regard, the notion of regularity can be seen as a kind of (stochastic) lower bound for $X_0$ when the system is supercritical. The notion of regularity, however, is more frequently used in Sections \ref{s:preliminaries}--\ref{s:2nd_crossing} and \ref{s:alpha}, where it is applied to critical systems. These critical systems either have a finite value of the corresponding $\Lambda(\cdot)$ (and they will be applied to prove Theorem \ref{t:main}), or are truncated critical systems whose values of $\Lambda(\cdot)$ are finite but depending on the level of truncation (and will be applied to prove Theorems \ref{t:2<alpha<4} and \ref{t:alpha=2}). In these cases, the notion of regularity should not be viewed as any kind of lower bound for (the initial distribution of) the system.\qed

\end{remark}

\medskip

Recall that under the integrability condition $\E(X_0 \, m^{X_0}) < \infty$, $p \ge p_c$ means $(m-1) \E(X_0 \, m^{X_0}) \ge  \E(m^{X_0})$.

\medskip

\begin{lemma}
 \label{A_4P(X_0=0)}

If $p\in [\, p_c, \, 1)$, then $(X_n, \, n\ge 0)$ is $0$-regular with coefficient $\chi=1-p$.
\end{lemma}

\medskip

\noindent {\it Proof.} Let $k\ge 1$ be an integer. By definition,
$$
\Xi_k(X_0)
=
\E \Big[ (X_0 \wedge k)^2 \, [(m-1)X_0-1] \, m^{X_0} \, {\bf 1}_{\{ X_0 \ge 1\} } \Big] .
$$

\noindent Since $(X_0 \wedge k)^2 \ge 1$ on $\{ X_0 \ge 1\}$, this yields
\begin{eqnarray*}
    \Xi_k(X_0)
 &\ge& \E \Big[ [(m-1)X_0-1] \, m^{X_0} \, {\bf 1}_{\{ X_0 \ge 1\} } \Big]
    \\
 &=& \E \Big[ [(m-1)X_0-1] \, m^{X_0} \Big] + \P (X_0=0) \, .
\end{eqnarray*}

\noindent Recall that $p\ge p_c$ means $\E \{ [(m-1)X_0-1] m^{X_0}\} \ge 0$, whereas $\P (X_0=0) = 1-p$. This implies $\Xi_k(X_0) \ge 1-p$, which is greater than or equal to $(1-p) \min\{ \Lambda (X_0), \, 1\}$. The lemma is proved.\qed

\medskip

\begin{remark}

 Let $\alpha \in [2, \, 4]$. Assume $\P(X_0^* = k) \sim c_0\, m^{-k} k^{-\alpha}$, $k\to \infty$, for some $0<c_0<\infty$. We are going to see in Lemma \ref{e:chi_lower_bound} that a critical system started at a conveniently truncated version of $X_0^*$ is $(4-\alpha)$-regular. This will allow to use a truncation argument in Section \ref{s:alpha} to prove Theorems \ref{t:2<alpha<4} and \ref{t:alpha=2}.\qed

\end{remark}

\section{Proof of Theorem \ref{t:main}: upper bound}
\label{s:ub}

The upper bound in Theorem \ref{t:main} does not require the full assumption in the theorem. In particular, here we do not need the assumption $\E(X_0^3 \, m^{X_0}) < \infty$.

Throughout the section, we assume $p>p_c>0$, i.e., $\E(m^{X_0}) < (m-1) \E(X_0 \, m^{X_0}) <\infty$. The upper bound in Theorem \ref{t:main} is as follows: there exists a constant $c_1>0$ such that for all sufficiently small $\varepsilon>0$,
\begin{equation}
    F_\infty (p_c+\varepsilon)
    \le
    \exp \Big( - \frac{c_1}{\varepsilon^{1/2}} \Big) .
    \label{ub}
\end{equation}

Let 
$$
\Lambda = \Lambda(X_0) := \E(X_0^3\, m^{X_0}) \in (0, \, \infty],
$$

\noindent as in \eqref{Lambda}. When $p\ge p_c$, we have $\E(m^{X_0}) \le (m-1) \E(X_0 \, m^{X_0})$ by definition, so $\Lambda \ge \E(X_0\, m^{X_0}) \ge \frac{1}{m-1} \E(m^{X_0}) \ge \frac{1}{m-1}$. Consequently,
\begin{equation}
    \Lambda \in [\frac{1}{m-1}, \, \infty] \quad \mathrm{if} \; p\ge p_c \, .
    \label{Lambda>}
\end{equation}

The main step in the proof of \eqref{ub} is the following proposition. We recall from Remark \ref{r:alpha=2;Lambda<infty} that if a system is $2$-regular in the sense of \eqref{alpha}, then $\Lambda<\infty$.

\medskip

\begin{proposition}
 \label{p:ub}

 Let $\beta\in [0, \, 2]$ and $\chi \in (0, \, 1]$.
 Assume $\E(X_0 \, m^{X_0}) <\infty$ and the system $(X_n, \, n\ge 0)$ is $\beta$-regular with coefficient $\chi$ in the sense of \eqref{alpha}. Let $\delta_0 := (m-1) \E(X_0 \, m^{X_0}) - \E(m^{X_0})$. There exist constants $c_2\in (0,\, 1]$ and $c_3\in (0,\, 1]$, depending only on $m$ and $\beta$, such that if $\delta_0 \in  (0, \, c_3)$ (for $\beta\in [0, \, 2)$) or if $\delta_0 \in  (0, \, c_3 \chi)$ (for $\beta=2$), then there exists an integer $n_0 \ge K_0$ satisfying
 \begin{equation}
      \E (m^{X_{n_0}}) \le 3,
      \qquad
      \max_{0\le i\le n_0} \E (X_i) \le 3 \, ,
      \label{eq_p:ub}
 \end{equation}

 \noindent where $K_0 := \min\{ (\frac{c_2 \, \chi \Lambda}{\delta_0})^{1/2}, \, (\frac{c_2\, \chi}{\delta_0})^{1/(2-\beta)}\}$ for $\beta\in [0, \, 2)$ and $K_0 := (\frac{c_2\chi\Lambda}{\delta_0})^{1/2}$ for $\beta=2$.

\end{proposition}

\medskip

By admitting Proposition \ref{p:ub} for the time being, we are able to prove the upper bound \eqref{ub} in Theorem \ref{t:main}.

\bigskip

\noindent {\it Proof of Theorem \ref{t:main}: upper bound.} Recall from \eqref{F_encadrement} that for all $n\ge 0$,
$$
    \frac{\E(X_n) - \frac{1}{m-1}}{m^n}
    \le
    F_\infty
    \le
    \frac{\E(X_n)}{m^n} \, ,
    \qquad
    n\ge 0\, .
$$

\noindent Lemma \ref{A_4P(X_0=0)} says that when $p\in (p_c, \, 1)$, the system $(X_n, \, n\ge 0)$ is $0$-regular with coefficient $\chi= 1-p$, so we are entitled to apply Proposition \ref{p:ub} to $\beta=0$. We take $p=p_c+ \varepsilon$ (with $0<\varepsilon< 1-p_c$); note that
\begin{eqnarray*}
    \delta_0
 &:=& (m-1) \E_{p_c+\varepsilon}(X_0 \, m^{X_0}) - \E_{p_c+\varepsilon}(m^{X_0})
    \\
 &=& [(m-1) \E_{p_c}(X_0 \, m^{X_0}) - \E_{p_c}(m^{X_0})]
    +
    \varepsilon\, (\E\{ [(m-1)X_0^*-1]m^{X_0^*}\} +1) \, .
\end{eqnarray*}

\noindent By Theorem A in the introduction, $(m-1) \E_{p_c}(X_0 \, m^{X_0}) - \E_{p_c}(m^{X_0}) =0$. So $\delta_0 = c_4 \, \varepsilon$ where $c_4 := \E\{ [(m-1)X_0^*-1]m^{X_0^*}\} +1 \in (1, \, \infty)$. By Proposition \ref{p:ub} (with $\beta=0$), $\E_{p_c+\varepsilon}(X_{\lfloor c_5/\varepsilon^{1/2}\rfloor}) \le 3$ for some constant $c_5>0$ and all sufficiently small $\varepsilon>0$. Using the second inequality in \eqref{F_encadrement}, we obtain the following bound for the free energy:  for all sufficiently small $\varepsilon>0$,
$$
F_\infty(p_c+\varepsilon)
\le
\frac{\E_{p_c+\varepsilon}(X_{\lfloor c_5/\varepsilon^{1/2}\rfloor})}{m^{\lfloor c_5/\varepsilon^{1/2}\rfloor}}
\le
3 \exp\Big( - (\frac{c_5}{\varepsilon^{1/2}} -1)\log m \Big),
$$

\noindent which yields \eqref{ub}.\qed

\bigskip

The rest of the section is devoted to the proof of Proposition \ref{p:ub}. Assume $\E(X_0^* \, m^{X_0^*}) <\infty$ (which is equivalent to saying that $p_c>0$) and $p>p_c$. For all $n\ge 0$, we write the moment generating function
$$
H_n(s)
:=
\E(s^{X_n}) \, .
$$

\noindent We rewrite the iteration equation \eqref{iteration} in terms of $H_n$: for all $n\ge 0$,
\begin{equation}
    H_{n+1}(s)
    =
    \frac1s \, H_n(s)^m
    +
    (1-\frac1s) \, H_n(0)^m \, .
    \label{Gn_iteration}
\end{equation}

\noindent A useful quantity in the proof is, for $n\ge 0$,
$$
\delta_n
:=
m(m-1) H_n'(m) - H_n(m)
=
(m-1) \E(X_n \, m^{X_n}) - \E(m^{X_n}) \, .
$$

\noindent By assumption, $\delta_0>0$ (and is small). Using the iteration relation \eqref{Gn_iteration}, it is immediate that
$$
\delta_{n+1}
=
H_n(m)^{m-1} \delta_n .
$$

\noindent Consequently, for $n\ge 1$,
\begin{equation}
    \delta_n
    =
    \Big( \prod_{i=0}^{n-1} H_i(m)^{m-1} \Big) \delta_0
    \in (0, \, \infty).
    \label{iteration_epsilon_n}
\end{equation}

\noindent [The recursion formula was known to Collet et al.~\cite{collet-eckmann-glaser-martin}; see also Equation (10) in \cite{bmxyz_questions}.] We outline the proof of Proposition \ref{p:ub} before getting into details.

\bigskip

\noindent {\it Outline of the proof of Proposition \ref{p:ub}.} Only the first inequality (saying that $\E(m^{X_{n_0}})\le 3$ for some integer $n_0 \ge K_0$) in the proposition needs to be proved. Let $\theta>0$ be a constant. We choose
$$
n_0 := \sup\{i\ge 1: \, \delta_i \le \theta\} \, .
$$

\noindent [If $n_0=\infty$, the proposition is easily proved. So we assume $n_0 <\infty$.] Since it is quite easy to see that $H_i(m) \le m^{1/(m-1)} \ee^{\delta_i}$ for $i\ge 0$, we get $H_{n_0}(m) \le m^{1/(m-1)} \ee^\theta\le 3$ if $\theta$ is chosen to satisfy $m^{1/(m-1)} \ee^\theta \le 3$. Consequently, $\E(m^{X_{n_0}})\le 3$ .

It remains to check that $n_0 \ge K_0$. The key ingredient is the following inequality: for all $0\le n\le n_0$,\footnote{Notation: $\prod_\varnothing := 1$.}
$$
\frac{\delta_n}{\delta_0}
=
\prod_{i=0}^{n-1} H_i(m)^{m-1} 
\le
\frac{c_6}{\Theta_0(s_n)} ,
\leqno(\ref{ub_pf_1})
$$

\noindent where $c_6>0$ is a constant depending only on $m$, $s_n := m(1- \frac{C_1}{n+1})$ for some constant $C_1\in (0, \, 1)$ depending on $m$, and $\Theta_i(\cdot)$, for $i\ge 0$, is a positive function on $(0, \, m)$ which is to be defined soon. [It turns out that for $s\in (0, \, m)$, $\Theta_0(s)$ is connected to $\Xi_k(X_0)$ defined in \eqref{Theta_0}, with $k = k(s):= \lfloor \frac{m}{m-s} \rfloor$. This helps explain partly the importance of the notion of regularity for systems.]

Let us look at \eqref{ub_pf_1} with $n:= n_0$. We have $\delta_{n_0} = \frac{\delta_{1+n_0}}{H_{n_0}(m)^{m-1}} > \frac{\theta}{H_{n_0}(m)^{m-1}}$. Since $H_{n_0}(m) \le 3$ as we have already seen, it follows that $\delta_{n_0}$ is greater than the positive constant $\frac{\theta}{3^{m-1}}$. Using a simple lower bound for $\Theta_0(s_n)$ (thanks to the aforementioned connection between $\Theta_0(s)$ and $\Xi_k(X_0)$) as a function of $n$, \eqref{ub_pf_1} will yield $n_0 \ge K_0$ as desired.

We are thus left with the proof of \eqref{ub_pf_1}. The idea is to study not only the function $\Theta_0(\cdot)$, but the sequence of functions $\Theta_n(\cdot)$, $n\ge 0$, and produce a recursion in $n$: for $0\le n \le n_0$, $s\in [\frac{m}{2}, \, m)$ and some constant $C_2>1$ depending only on $m$,
$$
\Theta_{n+1}(s) 
\ge
[C_2 s -(C_2-1)m] \, \frac{H_n(s)^{m-1}}{s} \, \Theta_n(s) \, .
$$

\noindent The proof of this inequality, done in two steps (Lemmas \ref{l:ub_2} and \ref{l:ub_3}), follows the lines of \cite{bmxyz_questions} and \cite{collet-eckmann-glaser-martin}. 

It is quite easy to see that $\frac{H_n(s)^{m-1}}{s} \ge (1-C_3(m-s))\frac{H_n(m)^{m-1}}{m}$ for some constant $C_3 >0$. So, as long as $s\in [\frac{m}{2}, \, m)$ satisfies $C_2 s -(C_2-1)m >0$ and $1-C_3(m-s)>0$ (which is the case when we choose $s:=s_n =m(1- \frac{C_1}{n+1})$ later), we get
$$
\Theta_{n+1}(s) 
\ge
\frac{[C_2 s -(C_2-1)m]\, (1-C_3(m-s))}{m} \, H_n(m)^{m-1} \, \Theta_n(s) \, ,
$$

\noindent and thus, by iteration,
$$
\Theta_n(s) 
\ge
\Big( \frac{[C_2 s -(C_2-1)m]\, (1-C_3(m-s))}{m}\Big)^n \, \Theta_0(s)\prod_{i=0}^{n-1} H_i(m)^{m-1} \, .
$$

\noindent On the other hand, $\Theta_n(s) \le 2$ for $0<s<m$ and $0\le n\le n_0$. With our choice of $s=s_n$, the factor $(\frac{[C_2 s -(C_2-1)m]\, (1-C_3(m-s))}{m})^n$ is greater than a positive constant independent of $n$. This will imply \eqref{ub_pf_1}.\qed

\bigskip

We now proceed to the detailed proof of Proposition \ref{p:ub}. Let us start with a few elementary properties of the moment generating functions.

\medskip

\begin{lemma}
 \label{l:ub_1}

 Assume $\E(X_0^* \, m^{X_0^*}) <\infty$ and $p>p_c$. Let $n\ge 0$.

 {\rm (i)} The functions $s\mapsto (m-1)s H_n'(s) - H_n(s)$ and $s\mapsto \frac{(m-1)s H_n'(s) - H_n(s)}{s}$ are non-decreasing on $(0, \, m]$.

 {\rm (ii)} We have $H_n(m) \le m^{1/(m-1)} \ee^{\delta_n}$.

 {\rm (iii)} For $s\in [1, \, m]$,
 $$
 \frac{H_n(s)^{m-1}}{s}
 \ge
 \frac{H_n(m)^{m-1}}{m} [1- m(m-s) \, \delta_n \, \ee^{(m-1)\delta_n}] \, .
 $$

\end{lemma}

\medskip

\noindent {\it Proof.} (i) Write
\begin{equation}
    \varphi_n(s)
    :=
    (m-1)sH_n'(s)-H_n(s) \, .
    \label{f_n}
\end{equation}

\noindent [So $\varphi_n(m) = \delta_n$.] For all $s\in (0, \, m)$,
$$
\varphi_n'(s)
=
(m-2) H_n'(s) + (m-1) s H_n''(s)
\ge
0 ,
$$

\noindent which implies the monotonicity of $\varphi_n$. In particular,
\begin{equation}
    \varphi_n(s)
    \le
    \delta_n,
    \qquad
    s\in (0, \, m)\, .
    \label{f_n<epsilon_n}
\end{equation}

To prove the monotonicity of $s\mapsto \frac{\varphi_n(s)}{s}$, we note that
$$
\frac{\!\d}{\! \d s} (\frac{\varphi_n(s)}{s})
=
\frac{s \varphi_n'(s) - \varphi_n(s)}{s^2}
=
\frac{(m-1)s^2 H_n''(s)+H_n(s) - sH_n'(s)}{s^2} .
$$

\noindent On the right-hand side, the numerator is greater than or equal to $s^2 H_n''(s)+H_n(s) - sH_n'(s)$, which is $\E \{ [X_n(X_n-1) +1 - X_n] s^{X_n}\} = \E \{ (X_n-1)^2 s^{X_n}\} \ge 0$, so the desired result follows.

(ii) For $s\in [1, \, m)$,
\begin{equation}
    \frac{\!\d}{\! \d s} (\frac{H_n(s)^{m-1}}{s})
    =
    [(m-1)s\, H_n'(s)-H_n(s)] \, \frac{H_n(s)^{m-2}}{s^2} \, .
    \label{pf:ub_eq0}
\end{equation}

\noindent By (i), $(m-1)s\, H_n'(s)-H_n(s) \le (m-1)m\, H_n'(m)-H_n(m) = \delta_n$. On the other hand, $\frac{H_n(s)^{m-2}}{s^2} \le \frac{H_n(s)^{m-1}}{s}$ (because $s\ge 1$ and $H_n(s) \ge 1$). Hence
\begin{equation}
    \frac{\!\d}{\! \d s} (\frac{H_n(s)^{m-1}}{s})
    \le
    \delta_n \, \frac{H_n(s)^{m-1}}{s}\, ,
    \qquad
    s\in [1, \, m)\, .
    \label{pf:ub_eq1}
\end{equation}

\noindent When $s=1$, $\frac{H_n(s)^{m-1}}{s}$ equals $1$, so solving the differential equation yields that for $s\in [1, \, m)$,
\begin{equation}
    \frac{H_n(s)^{m-1}}{s}
    \le
    \ee^{\delta_n (s-1)} \, .
    \label{pf:ub_eq2}
\end{equation}

\noindent Taking $s=m$ yields the desired inequality.

For further use, we also observe that our proof yields the following inequality in the {\it critical} case: if $(m-1) \E(Y_0 \, m^{Y_0}) = \E(m^{Y_0}) < \infty$, then
\begin{equation}
    \E(m^{Y_n}) \le m^{1/(m-1)},
    \qquad
    n\ge 0 \, .
    \label{Gn<}
\end{equation}

(iii) For $s\in [1,\, m)$, we have, by \eqref{pf:ub_eq1} and \eqref{pf:ub_eq2}, $\frac{\!\d}{\! \d s} (\frac{H_n(s)^{m-1}}{s}) \le \delta_n \ee^{(m-1)\delta_n}$. So for $s\in [1,\, m]$,
$$
\frac{H_n(s)^{m-1}}{s}
\ge
\frac{H_n(m)^{m-1}}{m}
-
(m-s) \, \delta_n \, \ee^{(m-1)\delta_n}\, .
$$

\noindent Since $H_n(m)^{m-1} \ge 1$, we have $(m-s) \, \delta_n \, \ee^{(m-1)\delta_n} \le \frac{H_n(m)^{m-1}}{m} m(m-s) \, \delta_n \, \ee^{(m-1)\delta_n}$, implying the desired conclusion.\qed

\bigskip

Define
$$
   \Theta_n(s)
   :=
   [H_n(s) - s(s-1)H_n'(s)]
   -
   \frac{(m-1)(m-s)}{m} \, [2sH_n'(s) + s^2 H_n''(s)]
   +
   \delta_n \, .
$$

\noindent [In the critical regime, $\delta_n=0$ for all $n$, so $\Theta_n(\cdot)$ is the function $\Delta_n(\cdot)$ studied in \cite{bmxyz_questions}.] Following the lines of \cite{bmxyz_questions} and \cite{collet-eckmann-glaser-martin}, we first prove two preliminary inequalities for $\Theta_n$.

\medskip

\begin{lemma}
 \label{l:ub_2}

 Assume $\E(X_0^* \, m^{X_0^*}) <\infty$, and $p>p_c$. Let $n\ge 0$.

 {\rm (i)} We have $\Theta_n(s)\in [0,\, 1+\delta_n]$ for all $s\in [0,\, m)$.

 {\rm (ii)} If $\delta_n \le 1$, then for all $s\in [\frac{m}{2}, \, m)$,
 $$
 \Theta_{n+1}(s)
    \ge
    \frac{m}{s}\, \Theta_n(s) \, H_n(s)^{m-1}
    -
    \frac{\kappa (m-s)}{s}\, [\delta_n - \varphi_n(s)]^2  H_n(s)^{m-2} ,
 $$
 where $\kappa := 3^{m-2}+1$, and $\varphi_n(s) := (m-1)sH_n'(s)-H_n(s)$ as in \eqref{f_n}.

\end{lemma}

\medskip

\noindent {\it Proof.} (i) By definition, for $s\in (0, \, m)$,
$$
\Theta_n'(s)
=
- \frac{m-s}{m} \, [ 2(m-2)H_n'(s) + (4m-5) sH_n''(s) + (m-1) s^2H_n'''(s)]
\le 
0\, ,
$$

\noindent so $s\mapsto \Theta_n(s)$ is non-increasing on $(0, \, m)$. Since $\lim_{s\to m-} \Theta_n(s) = 0$, and $\Theta_n(0) = H_n(0)+\delta_n \le 1+ \delta_n$, the result follows.

(ii) Let $s\in [\frac{m}{2}, \, m)$. The iteration \eqref{Gn_iteration} yields
\begin{eqnarray}
    \Theta_{n+1}(s)
 &=& \frac{m}{s}\, \Theta_n(s) \, H_n(s)^{m-1}
    -
    \frac{m-s}{s}\, \varphi_n(s)^2 H_n(s)^{m-2}
    \nonumber
    \\
 && \qquad
    +
    m(\frac{H_n(m)^{m-1}}{m}-\frac{H_n(s)^{m-1}}{s})\delta_n .
    \label{pf:ub_eq3}
\end{eqnarray}

\noindent For the last term on the right-hand side, we recall from \eqref{pf:ub_eq0} that $\frac{\!\d}{\! \d s} (\frac{H_n(s)^{m-1}}{s}) = \varphi_n(s) \, \frac{H_n(s)^{m-2}}{s^2}$, so by the mean-value theorem, there exists $y\in [s, \, m)$ such that
$$
\frac{H_n(m)^{m-1}}{m}-\frac{H_n(s)^{m-1}}{s}
=
(m-s) \, \frac{\varphi_n(y)}{y} \, \frac{H_n(y)^{m-2}}{y}
\ge
(m-s) \, \frac{\varphi_n(s)}{s}\, \frac{H_n(y)^{m-2}}{y} \, ,
$$

\noindent the last inequality being a consequence of Lemma \ref{l:ub_1}~(i). Going back to \eqref{pf:ub_eq3}, we see that the proof will be finished if we are able to check that for all $u\in [s, \, m)$, with $\kappa := 3^{m-2}+1$,
\begin{equation}
    \varphi_n(s)^2  H_n(s)^{m-2}
    -
    \frac{m}{u} \, \varphi_n(s) \, H_n(u)^{m-2}\, \delta_n
    \le
    \kappa [\delta_n - \varphi_n(s)]^2  H_n(s)^{m-2} \, .
    \label{pf:ub_eq4}
\end{equation}

We prove \eqref{pf:ub_eq4} by distinguishing two possible situations. We write $\mathrm{LHS}_{\eqref{pf:ub_eq4}}$ for the expression on the left-hand side of \eqref{pf:ub_eq4}, and $\mathrm{RHS}_{\eqref{pf:ub_eq4}}$ for the expression on the right-hand side.

First situation: $\varphi_n(s) \ge 0$. We use the trivial inequalities $\frac{m}{u} \ge 1$ and $H_n(u)^{m-2} \ge H_n(s)^{m-2}$, to see that
$$
\mathrm{LHS}_{\eqref{pf:ub_eq4}}
\le
\varphi_n(s)^2 H_n(s)^{m-2}
-
\varphi_n(s) H_n(s)^{m-2} \delta_n \, ,
$$

\noindent which is non-positive (because $\varphi_n(s) \le \delta_n$; see \eqref{f_n<epsilon_n}). This yields \eqref{pf:ub_eq4} since $\mathrm{RHS}_{\eqref{pf:ub_eq4}} \ge 0$.

Second (and last) situation: $\varphi_n(s) <0$. We write $|\varphi_n(s)|$ instead of $-\varphi_n(s)$ in this situation. We have
$$
H_n(m)
\le
H_n(s) + (m-s) H_n'(m)
=
H_n(s) + (m-s) \frac{H_n(m) + \delta_n}{m(m-1)} \, .
$$

\noindent For $s\in [\frac{m}{2}, \, m)$, we have $\frac{m-s}{m(m-1)} \le \frac{m-s}{m} \le \frac12$, so
$$
H_n(m)
\le
H_n(s) + \frac{H_n(m) + \delta_n}{2} \, .
$$

\noindent Consequently,
\begin{equation}
    H_n(m)
    \le
    2 H_n(s) + \delta_n
    \le
    2 H_n(s) + 1
    \le 3 H_n(s) \, ,
    \label{pf:ub_eq4-bis}
\end{equation}

\noindent where we used the assumption $\delta_n \le 1$ in the second inequality, and the trivial relation $H_n(s) \ge 1$ in the last inequality.

We now look at the second expression on the left-hand side of \eqref{pf:ub_eq4}. The factor $\frac{m}{u}$ is easy to deal with: we have $\frac{m}{u} \le 2$ (using $u \ge s \ge \frac{m}{2}$). For the factor $H_n(u)^{m-2}$, since $u\le m$, we have $H_n(u) \le H_n(m) \le 3 H_n(s)$ (by means of \eqref{pf:ub_eq4-bis}), so $H_n(u)^{m-2} \le 3^{m-2} H_n(s)^{m-2}$. Consequently,
$$
\mathrm{LHS}_{\eqref{pf:ub_eq4}}
\le
|\varphi_n(s)|^2  H_n(s)^{m-2}
+
2|\varphi_n(s)| \, 3^{m-2} H_n(s)^{m-2} \, \delta_n\, .
$$

\noindent We look at the two terms on the right-hand side. For the first term, we argue that $|\varphi_n(s)|^2 \le (|\varphi_n(s)|+ \delta_n)^2$. For the second term, we use $2|\varphi_n(s)| \, \delta_n \le (|\varphi_n(s)|+ \delta_n)^2$. Hence
\begin{eqnarray*}
    \mathrm{LHS}_{\eqref{pf:ub_eq4}}
 &\le& (|\varphi_n(s)|+ \delta_n)^2 H_n(s)^{m-2}
    +
    (|\varphi_n(s)|+ \delta_n)^2 \, 3^{m-2} H_n(s)^{m-2}
    \\
 &=& (3^{m-2}+1) (|\varphi_n(s)|+ \delta_n)^2 H_n(s)^{m-2} \, ,
\end{eqnarray*}

\noindent which yields \eqref{pf:ub_eq4} again, as $[\delta_n -\varphi_n(s)]^2 = (|\varphi_n(s)|+ \delta_n)^2$ in this case. [Note that this case is very easy to handle when $m=2$: all we need is to observe that $\frac{m}{u} \le 2$.]\qed

\medskip

\begin{lemma}
 \label{l:ub_3}

 Assume $\E(X_0 \, m^{X_0}) <\infty$, and $p>p_c$. Let $n\ge 0$. Then
 $$
 [\delta_n - \varphi_n(s)]^2
 \le
 2(H_n(0)+\delta_n) \Theta_n(s),
 \qquad
 s\in (0, \, m) ,
 $$
 where $\varphi_n(s) := (m-1)sH_n'(s)-H_n(s)$ as in \eqref{f_n}.

\end{lemma}

\medskip

\noindent {\it Proof.} The lemma in the critical regime was already proved in (\cite{bmxyz_questions}, proof of Lemma 7). The argument remains valid in our situation if we replace $G_n(s)$ and $\Delta_n(s)$ there (notation of \cite{bmxyz_questions}) by $H_n(s)+\delta_n$ and $\Theta_n(s)$, respectively. It is reproduced here for the sakes of clarity and self-containedness.

By definition, with $t:= \frac{s}{m} \in (0, \, 1)$,
\begin{eqnarray}
    \Theta_n(s)
 &=&[H_n(s) - s(s-1)H_n'(s)] - [H_n(m) - m(m-1)H_n'(m)]
    \nonumber
    \\
 &&\qquad
    -
    \frac{(m-1)(m-s)}{m} \, [2sH_n'(s) + s^2 H_n''(s)]
    \nonumber
    \\
 &=& \E \Big\{ m^{X_n} [(m-1)X_n-1] [1-(1+X_n)t^{X_n} + X_n \, t^{1+X_n}] \Big\}
    \nonumber
    \\
 &=& \E \Big\{ m^{X_n} [(m-1)X_n-1] [1-(1+X_n)t^{X_n} + X_n \, t^{1+X_n}]\, {\bf 1}_{\{ X_n \ge 1\} } \Big\} ,
    \label{Theta_n}
\end{eqnarray}

\noindent and
\begin{eqnarray*}
    \delta_n - \varphi_n(s)
 &=& \E \Big\{ m^{X_n} (1-t^{X_n}) [(m-1)X_n-1]\, {\bf 1}_{\{ X_n \ge 1\} } \Big\} ,
    \\
    H_n(0)+\delta_n
 &=& \E \Big\{ m^{X_n} [(m-1)X_n-1]\, {\bf 1}_{\{ X_n \ge 1\} } \Big\} \, .
\end{eqnarray*}

\noindent By the Cauchy--Schwarz inequality,
$$
[\delta_n - \varphi_n(s)]^2
\le
(H_n(0)+\delta_n) \E \Big\{ m^{X_n} (1-t^{X_n})^2 [(m-1)X_n-1] \, {\bf 1}_{\{ X_n \ge 1\} } \Big\}.
$$

\noindent So the proof of the lemma is reduced to showing the following: for $t \in (0, \, 1)$ and integer $k\ge 1$,
$$
(1-t^k)^2 \le 2[1-(1+k)t^k+kt^{1+k}] \, .
$$

\noindent This is equivalent to saying that $2 k\, t^k(1-t)\le 1-t^{2k}$, which is obviously true because $\frac{1-t^{2k}}{1-t}= 1+t+...+t^{2k-1} \ge 2k t^k$ (noting that $2t^k\le 2t^{k-1/2} \le t^{k+\ell}+ t^{k-\ell-1}$ for $0\le \ell < k$).\qed

\bigskip

We have now all the ingredients for the proof of Proposition \ref{p:ub}.

\bigskip

\noindent {\it Proof of Proposition \ref{p:ub}.} Assume $\delta_0:=(m-1)\E(X_0m^{X_0})-\E(m^{X_0})>0$. If we are able to prove the first inequality (saying that $\E (m^{X_{n_0}}) \le 3$ for some integer $n_0 \ge K_0$), then $\E (X_{n_0}) \le \E (m^{X_{n_0}}) \le 3$. Since $m \, \E(X_n) -1 \le \E (X_{n+1})$, i.e., $\E (X_n) \le \frac{1+ \E (X_{n+1})}{m}$ (for all $n\ge 0$), the second inequality in Proposition \ref{p:ub} follows immediately.

It remains to prove the first inequality. Fix a constant $0<\theta \le 1$ such that $m^{1/(m-1)} \ee^\theta\le 3$. Let us define
$$
n_0
=
n_0 (\theta, \, \delta_0)
:=
\sup \{ i\ge 0: \, \delta_i \le \theta \} \, .
$$

\noindent We may assume that $\delta_0\in (0, \, \frac{\theta}{m\ee})$, so by Lemma \ref{l:ub_1}~(ii), $\delta_1 = H_0(m)^{m-1} \delta_0 \le[m^{1/(m-1)}\ee^{\delta_0} ]^{(m-1)} \delta_0 < \theta$. This implies $n_0 \ge 1$. 

If $n_0=\infty$, then $\delta_i \le \theta$ for all $i\ge 0$, and by Lemma \ref{l:ub_1}~(ii), $\E(m^{X_n})\le m^{1/(m-1)} \ee^{\delta_n}\le m^{1/(m-1)} \ee^\theta\le 3$ for all $n\ge 0$: there is nothing to prove in the proposition. In the following, we assume that $n_0<\infty$. 

Let $0\le n \le n_0$ (hence $\delta_n \le \theta \le 1$) and $s\in (\frac{4\kappa-1}{4\kappa}m, \, m)$, where $\kappa := 3^{m-2}+1$ as in Lemma \ref{l:ub_2}. We reproduce the conclusion of Lemma \ref{l:ub_2}~(ii) in its notation:
$$
\Theta_{n+1}(s)
\ge
\frac{m}{s}\, \Theta_n(s) \, H_n(s)^{m-1}
-
\frac{\kappa (m-s)}{s}\, [\delta_n - \varphi_n(s)]^2  H_n(s)^{m-2} .
$$

\noindent By Lemma \ref{l:ub_3}, $[\delta_n - \varphi_n(s)]^2 \le 2(H_n(0)+\delta_n) \Theta_n(s)$; since $H_n(0) \le 1$ and $\delta_n \le \theta \le 1$ (for $0\le n\le n_0$), we have $[\delta_n - \varphi_n(s)]^2 \le 4\Theta_n(s) \le 4 \Theta_n(s) H_n(s)$ (because $H_n(s) \ge 1$ for $s\ge 1$); this implies that for $0\le n\le n_0$,
\begin{eqnarray*}
    \Theta_{n+1}(s)
 &\ge& \frac{m}{s}\, \Theta_n(s) \, H_n(s)^{m-1}
    -
    \frac{4\kappa(m-s)}{s}\, \Theta_n(s)  H_n(s)^{m-1}
    \\
 &=& \frac{4\kappa s-(4\kappa-1)m}{s} \, H_n(s)^{m-1}\Theta_n(s).
\end{eqnarray*}

\noindent Iterating the procedure, we see that for $0\le n\le n_0$ and $s\in (\frac{4\kappa-1}{4\kappa}m, \, m)$,
$$
\Theta_n (s)
\ge
\Theta_0(s) \, [4\kappa s-(4\kappa-1)m]^n \prod_{i=0}^{n-1} \frac{H_i(s)^{m-1}}{s} \, .
$$

\noindent By Lemma \ref{l:ub_1}~(iii),
$$
\frac{H_n(s)^{m-1}}{s}
\ge
\frac{H_n(m)^{m-1}}{m} [1- m(m-s) \, \delta_n \, \ee^{(m-1)\delta_n}] \, ,
$$

\noindent which is greater than or equal to $\frac{H_n(m)^{m-1}}{m} [1- m(m-s) \, \theta\, \ee^{(m-1)\theta}]$ for $0\le n\le n_0$. Hence, for $0\le n\le n_0$ and $s\in (\frac{4\kappa-1}{4\kappa}m, \, m)$ satisfying $m(m-s) \, \theta \, \ee^{(m-1)\theta}<1$,
$$
\prod_{i=0}^{n-1} H_i(m)^{m-1}
\le
\frac{\Theta_n (s)}{\Theta_0 (s)} \,
\frac{m^n}{[4\kappa s-(4\kappa-1)m]^n} \,
\frac{1}{[1- m(m-s) \, \theta \, \ee^{(m-1)\theta}]^n} \, .
$$

\noindent We take $s = s_n:= (1-\frac{1}{8m^2\kappa (n+1)})m \in (\frac{4\kappa-1}{4\kappa}m, \, m)$ from now on. The requirement $m(m-s_n) \, \theta \, \ee^{(m-1)\theta}<1$ is met with our choice: since $\theta\le 1$ and $\ee^\theta \le 3$, we have
$$
m(m-s_n) \, \theta \, \ee^{(m-1)\theta}
=
\frac{\theta \, \ee^{(m-1)\theta}}{8\kappa(n+1)}
\le
\frac{3^{m-1}}{8(3^{m-2}+1)(n+1)}
<
1.
$$

\noindent By Lemma \ref{l:ub_2}~(i), $\Theta_n (s_n) \le 1+ \delta_n \le 1+ \theta \le 2$ (for $0\le n\le n_0$). This yields the existence of a constant $c_6>0$, depending only on $m$, such that for $0\le n\le n_0$,
\begin{equation}
    \prod_{i=0}^{n-1} H_i(m)^{m-1}
    \le
    \frac{c_6}{\Theta_0 (s_n)} \, .
    \label{ub_pf_1}
\end{equation}

Let us have a closer look at $\Theta_0(s_n)$. Recall from \eqref{Theta_n} that for $s\in (1, \, m)$ and $t:= \frac{s}{m} \in (0, \, 1)$,
$$
\Theta_0(s)
=
\E \Big\{ m^{X_0} [(m-1)X_0-1] [1-(1+X_0)t^{X_0} + X_0 \, t^{1+X_0}]\, {\bf 1}_{\{ X_0 \ge 1\} } \Big\} .
$$

\noindent For $t\in (0, \, 1)$ and $k\ge 1$, we have $t^k \le \ee^{-(1-t)k}$, so $1-(1+k)t^k+ kt^{1+k} \ge 1-(1+u)\, \ee^{-u}$, where $u:= (1-t)k>0$. Observe that $1-(1+v)\, \ee^{-v} \ge 1- \frac{2}{\ee}$ for $v\ge 1$ (because $v\mapsto 1-(1+v)\, \ee^{-v}$ is increasing on $(0, \, \infty)$) and that $1-(1+v)\, \ee^{-v} \ge \frac{v^2}{2\ee}$ for $v\in (0, \, 1]$ (because $v\mapsto 1-(1+v)\, \ee^{-v}-\frac{v^2}{2\ee}$ is increasing on $(0, \, 1]$). Hence
$$
1-(1+k)t^k+ kt^{1+k}
\ge
c_7 \, \min\{ (1-t)^2k^2, \, 1\} ,
\qquad
t\in (0, \, 1), \; k\ge 1 \, ,
$$

\noindent with $c_7 := \min\{ 1- \frac{2}{\ee}, \, \frac{1}{2\ee}\} >0$. Consequently, for all $s\in (0, \, m)$,
\begin{eqnarray*}
    \Theta_0(s)
 &\ge& c_7 \, \E \Big\{ m^{X_0} [(m-1)X_0-1] \min\{ (1-\frac{s}{m})^2X_0^2, \, 1\}\, {\bf 1}_{\{ X_0 \ge 1\} } \Big\}
    \\
 &=& c_7 \, (1-\frac{s}{m})^2 \, \E \Big\{ m^{X_0} [(m-1)X_0-1] \min\{ X_0^2, \, \frac{m^2}{(m-s)^2}\} \Big\} .
\end{eqnarray*}

\noindent Let $\Xi_k(X_0) := \E \{ (X_0 \wedge k)^2 \, [(m-1)X_0-1] \, m^{X_0} \}$ be as in \eqref{Theta_0}. Then this gives
$$
\Theta_0(s)
\ge
c_7 \, (1-\frac{s}{m})^2 \, \Xi_{\lfloor \frac{m}{m-s}\rfloor} (X_0) \, .
$$

\noindent We now come back to our choice of $s=s_n:= (1-\frac{1}{8m^2\kappa (n+1)})m$, so $1-\frac{s}{m} = \frac{1}{8m^2\kappa (n+1)}$, and $\lfloor \frac{m}{m-s}\rfloor = \frac{m}{m-s} = 8m^2\kappa (n+1)$. By assumption, the system $(X_n, \, n\ge 0)$ is $\beta$-regular with coefficient $\chi$ in the sense of \eqref{alpha}, i.e., $\Xi_k(X_0) \ge \chi\, \min\{ \Lambda, \, k^\beta\}$ for all integers $k\ge 1$. Consequently, for $n\ge 0$, with $c_8 := \frac{c_7}{(8m^2 \kappa)^2}$,
$$
\Theta_0(s_n)
\ge
\frac{c_8 \, \chi}{(n+1)^2}\, \min\{ \Lambda, \, (n+1)^\beta\} \, .
$$

\noindent Combining with \eqref{ub_pf_1}, we get, for $0\le n\le n_0$ and with $c_9 := \frac{c_6}{c_8}$,
$$
\prod_{i=0}^{n-1} H_i(m)^{m-1}
\le
\frac{c_9}{\chi}\, \frac{(n+1)^2}{\min\{ \Lambda, \, (n+1)^\beta\}} \, .
$$

\noindent Recall from \eqref{iteration_epsilon_n} that $\delta_n = (\prod_{i=0}^{n-1} H_i(m)^{m-1} ) \delta_0$, which yields: for $0\le n\le n_0$,
$$
\delta_n
\le
\frac{c_9}{\chi}\, \frac{(n+1)^2\, \delta_0}{\min\{ \Lambda, \, (n+1)^\beta\}} \, .
$$

\noindent We take $n=n_0$. Then on the one hand, we have
\begin{equation}
    \delta_{n_0}
    \le
    \frac{c_9}{\chi}\, \frac{(n_0+1)^2\, \delta_0}{\min\{ \Lambda, \, (n_0+1)^\beta\}} \, .
    \label{pf:ub_eq6}
\end{equation}

\noindent On the other hand,
$$
\delta_{n_0}
=
\frac{\delta_{n_0+1}}{H_{n_0}(m)^{m-1}}
\ge
\frac{\theta}{H_{n_0}(m)^{m-1}} \, .
$$

\noindent By Lemma \ref{l:ub_1}~(ii), $H_{n_0}(m) \le m^{1/(m-1)} \ee^{\delta_{n_0}} \le m^{1/(m-1)} \ee^\theta$, which is bounded by $3$ by the choice of $\theta$; hence
$$
\E(m^{X_{n_0}}) \le 3 \, .
$$

\noindent We have $\delta_{n_0} \ge \frac{\theta}{H_{n_0}(m)^{m-1}} \ge \frac{\theta}{m \, \ee^{(m-1)\theta}}$. Combining this with \eqref{pf:ub_eq6}, we get
$$
\frac{(n_0+1)^2}{\min\{ \Lambda, \, (n_0+1)^\beta\}}
\ge
\frac{c_{10} \, \chi}{\delta_0}\, ,
$$

\noindent with $c_{10} := \frac{\theta}{c_9 \, m \, \ee^{(m-1)\theta}}$. The proposition is proved (recalling from Remark \ref{r:alpha=2;Lambda<infty} that $\Lambda<\infty$ in case $\beta=2$).\qed

\section{Open paths, $XY$-coupling, and the first crossing}
\label{s:preliminaries}

The proof of the lower bound in Theorems \ref{t:main}, \ref{t:2<alpha<4} and \ref{t:alpha=2} requires deep probabilistic understanding of the model in the supercritical regime. It is made possible via a coupling argument, called $XY$-coupling, by means of a comparison between the supercritical system and an appropriate critical system. The coupling is quite elementary if we use the hierarchical representation of recursive systems.

The organization of this section is as follows. Section \ref{subs:hierarchical} introduces the hierarchical representation of recursive systems, regardless of whether they are supercritical or critical (or even subcritical, though we do not deal with subcritical systems in this paper). Section \ref{subs:coupling} gives the $XY$-coupling between the supercritical system and an appropriate critical system. In the brief but crucial Section \ref{subs:open_paths}, we introduce an important quantity for the critical system: the number of {\it open paths}. Section \ref{subs:1st-key} is a bridge connecting the expected value in a supercritical system and the expected number of open paths in a critical system. This bridge is crossed a first time in Section \ref{subs:the_first_crossing} to study a critical system by means of an appropriate supercritical system. [In Section \ref{s:2nd_crossing}, this bridge will be crossed a second time, in the opposite direction, to study a supercritical system by means of an appropriate critical system. Both crossings need much preparation: the first crossing relies on technical notions (hierarchical representation, open paths and the $XY$-coupling), whereas the preparation for the second crossing involves some technical estimates in Sections \ref{s:open_paths} and \ref{s:2nd_crossing}.]

\subsection{The hierarchical representation}
\label{subs:hierarchical}

Let $(X_n, \, n\ge 0)$ be a system as defined in \eqref{iteration}, i.e., such that each $X_{n+1}$ has the law of $(X_{n,1} + \cdots + X_{n,m} -1)^+$, where $X_{n,k}$, $1\le k\le m$, are independent copies of $X_n$. In order to introduce the $XY$-coupling between a supercritical system and a critical system, it turns out to be convenient to use a simple hierarchical representation of the system. As a matter of fact, it is in the form of the hierarchical representation that the system appeared in Collet et al.~\cite{collet-eckmann-glaser-martin2} and in Derrida and Retaux~\cite{derrida-retaux}.

We define a family of random variables $(X(x), \, x\in \T)$, indexed by a (reversed) $m$-ary tree $\T$, in the following way. For any vertex $x$ in the genealogical tree $\T$, we use $|x|$ to denote the generation of $x$; so $|x|=0$ if $x$ is in the initial generation. We assume that $X(x)$, for $x\in \T$ with $|x|=0$ (i.e., in the initial generation of $\T$), are i.i.d.\ having the distribution of $X_0$. For any vertex $x\in \T$ with $|x|\ge 1$, let $x^{(1)}$, $\ldots$, $x^{(m)}$ denote the $m$ parents of $x$ in generation $|x|-1$, and set
$$
X(x)
:=
(X(x^{(1)}) + \cdots + X(x^{(m)})-1)^+ \, .
$$

\begin{figure}[htb]
\centering{\includegraphics[width=.75\textwidth]{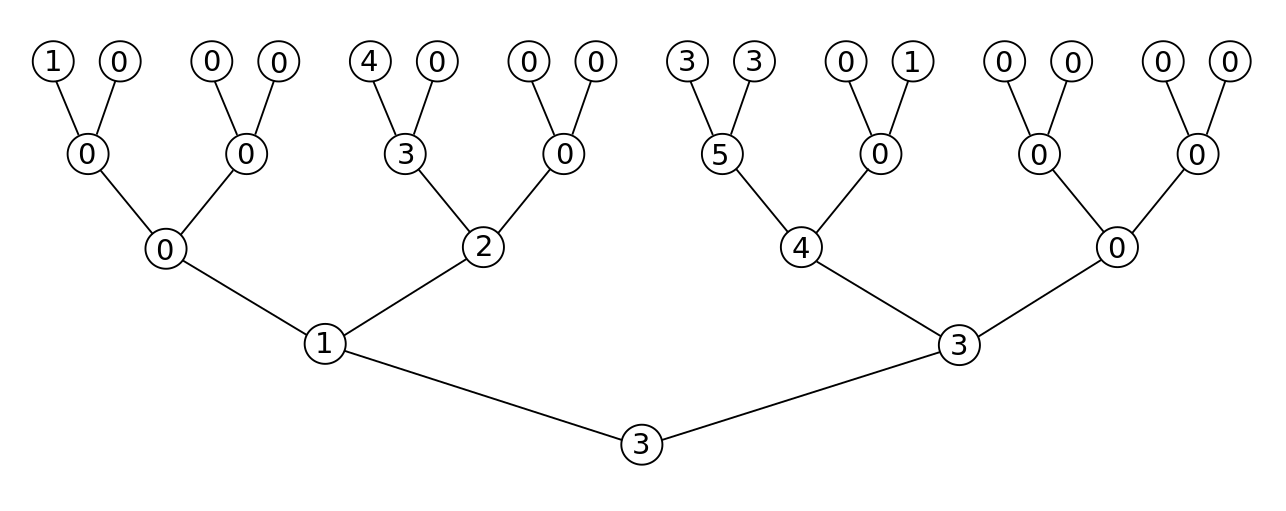}}
\caption{\leftskip=1.8truecm \rightskip=1.8truecm An example of hierarchical representation (the first 4 generations) with $m=2$. Each vertex is represented by a circle; the number in the circle is the value of the system at the vertex.}

 \label{f:fig1}
\end{figure}

\noindent See Figure \ref{f:fig1}.

As such, for any given $n\ge 0$, the sequence of random variables $X(x)$, for $x\in \T$ with $|x|=n$, are i.i.d.\ having the distribution of $X_n$.

The hierarchical representation is valid for any system satisfying the recursion \eqref{iteration}, regardless of whether the system is supercritical, critical or subcritical.

\subsection{The $XY$-coupling}
\label{subs:coupling}

Let $(X_n, \, n\ge 0)$ be a system as defined in \eqref{iteration}. We are going to make a coupling, called $XY$-coupling, for $(X_n, \, n\ge 0)$ and an appropriate system $(Y_n, \, n\ge 0)$ such that the first system is always greater than or equal to the second in their hierarchical representation.

Assume $\P(X_0=k) \ge \P(Y_0 =k)$ for all integers $k\ge 1$, and $\P(Y_0=0)>0$. We can couple the random variables $X_0$ and $Y_0$ in a same probability space such that $X_0 \ge Y_0$ a.s.\ and that $\P(X_0 = Y_0 \, | \, Y_0>0) =1$.\footnote{Let $Z$ be a $\z_+$-valued random variable, independent of $Y_0$, such that $\P(Z=k) := \frac{\P(X_0=k)-\P(Y_0=k)}{\P(Y_0=0)}$ for $k\ge 1$ and $\P(Z=0) := \frac{\P(X_0=0)}{\P(Y_0=0)}$. Then the pair $Y_0 + Z\, {\bf 1}_{\{ Y_0=0\} }$ and $Y_0$ will do the job.} 

For each $n\ge 0$, let $Y_{n+1}$ be a random variable having the law of $(Y_{n,1} + \cdots + Y_{n,m} -1)^+$, where $Y_{n,1}$, $\ldots$, $Y_{n,m}$ are independent copies of $Y_n$.

As in Section \ref{subs:hierarchical}, let $\T$ be the common genealogical $m$-ary tree associated with systems $(X_n, \, n\ge 0)$ and $(Y_n, \, n\ge 0)$. Let $(X(x), \, Y(x))$, for $x\in \T$ with $|x|=0$, be i.i.d.\ random variables having the distribution of $(X_0, \, Y_0)$, where $X_0$ and $Y_0$ are already coupled such that $X_0\ge Y_0$ a.s.\ and that $\P(X_0=Y_0 \, | \, Y_0 >0) =1$. For any vertex $x\in \T$ with $|x|\ge 1$, let
\begin{eqnarray*}
    X(x)
 &:=& (X(x^{(1)}) + \cdots + X(x^{(m)})-1)^+ ,
    \\
    Y(x)
 &:=& (Y(x^{(1)}) + \cdots + Y(x^{(m)})-1)^+ \, ,
\end{eqnarray*}

\noindent where $x^{(1)}$, $\ldots$, $x^{(m)}$ denote as before the parents of $x$. It follows that for any given $n\ge 0$, the sequence of random variables $X(x)$, indexed by $x\in \T$ with $|x|=n$, are i.i.d.\ having the distribution of $X_n$, whereas the sequence of random variables $Y(x)$, also indexed by $x\in \T$ with $|x|=n$, are i.i.d.\ having the distribution of $Y_n$. Moreover, for all $x\in \T$, $X(x) \ge Y(x)$ a.s.

We observe that since $X_0 \ge Y_0$ a.s.\ and $\P( X_0 = Y_0 \, | \, Y_0 >0)=1$, we have
$$
\E(X_0-Y_0 \, |\, Y_0=0)
=
\frac{\E[(X_0-Y_0)\, {\bf 1}_{\{ Y_0=0\} }]}{\P(Y_0=0)}
=
\frac{\E(X_0-Y_0)}{\P(Y_0=0)} \, .
$$

Here is a summary of the $XY$-coupling.

\medskip

\begin{theorem}
\label{t:coupling}

 {\bf (The $XY$-coupling).} Assume $\P(X_0=k) \ge \P(Y_0 =k)$ for all integers $k\ge 1$, and $\P(Y_0=0)>0$. We can couple two systems $(X(x), \, x\in \T)$ and $(Y(x), \, x\in \T)$ with initial distributions $X_0$ and $Y_0$, respectively, such that $X(x) \ge Y(x)$ for all $x\in \T$; in particular, we can couple two systems $(X_n, \, n\ge 0)$ and $(Y_n, \, n\ge 0)$ such that $X_n \ge Y_n$ \hbox{\rm a.s.} for all $n\ge 0$. Moreover, the coupling satisfies
 \begin{equation}
     \E(X_0-Y_0 \, |\, Y_0=0)
     =
     \frac{\E(X_0-Y_0)}{\P(Y_0=0)} \, .
     \label{E(Z|Y)}
 \end{equation}

\end{theorem}

\medskip

We are going to apply the $XY$-coupling several times. Each time, $(Y_n, \, n\ge 0)$ is critical (so the condition $\P(Y_0=0)>0$ is automatically satisfied), and $(X_n, \, n\ge 0)$ supercritical.

\subsection{Open paths}
 \label{subs:open_paths}

Let $(Y_n, \, n\ge 0)$ denote a system satisfying \eqref{iteration}. For any vertex $x\in \T$, we call $(x_k, \, 0\le k\le |x|)$ a path leading to $x$ if each $x_{k+1}$ is the (unique) child of $x_k$, and $|x_k| = k$. [Degenerate case: when $|x|=0$, the path leading to $x$ is reduced to the singleton $x$.] A path is said to be {\it open} if for any vertex $x$ on the path with $|x|\ge 1$, we have $Y(x) = Y(x^{(1)}) + \cdots + Y(x^{(m)})-1$ (or equivalently, $Y(x^{(1)}) + \cdots + Y(x^{(m)}) \ge 1$). [Degenerate case: when $|x|=0$, the path is considered as open.]

For $x\in \T$ and integer $i\ge 0$, let $N^\# (x)$ denote the number of open paths leading to $x$, and $N^{(i)}(x)$ the number of open paths leading to $x$ such that $Y(x_0)=i$; note that an open path can start with a vertex $x_0$ in the initial generation with $Y(x_0)=0$ if it receives enough support from neighbouring vertices along the path. By definition, $0\le N^{(i)} (x) \le N^\# (x) \le m^{|x|}$ for all $x\in \T$ and integers $i\ge 0$, and $N^\# (y) = 1$ and $N^{(i)}(y) = {\bf 1}_{\{ Y(y)=i\} }$ for $y\in \T$ with $|y|=0$.

\begin{figure}[htb]
\centering \setlength\unitlength{0.75pt}
\begin{picture}(400,160)
\put(450,10)\centering{\includegraphics[width=.85\textwidth]{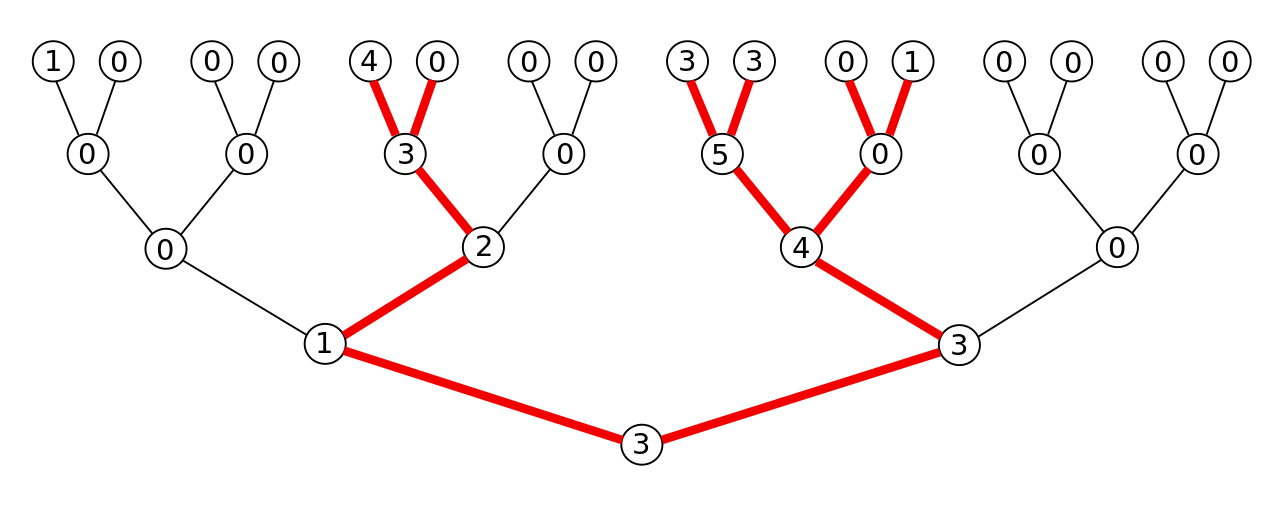}}
\put(-414,140){\small ${\mathfrak e}_0$}
\put(-402,110){\small ${\mathfrak e}_1$}
\put(-378,81){\small ${\mathfrak e}_2$}
\put(-328,51){\small ${\mathfrak e}_3$}
\put(-228,18){\small ${\mathfrak e}_4$}
\end{picture}\caption{\leftskip=1.6truecm \rightskip=1.6truecm The same example as in Figure \ref{f:fig1}, with open paths marked in bold (and coloured in red). For this example, $N_4^{\#} = 6$, $N_4^{(0)} = 2$, $N_4^{(1)} =1$, $N_4^{(2)} =0$, $N_4^{(3)}=2$, $N_4^{(4)} =1$, and $N_4^{(i)}=0$ if $i\ge 5$.}

 \label{f:fig2}
\end{figure}
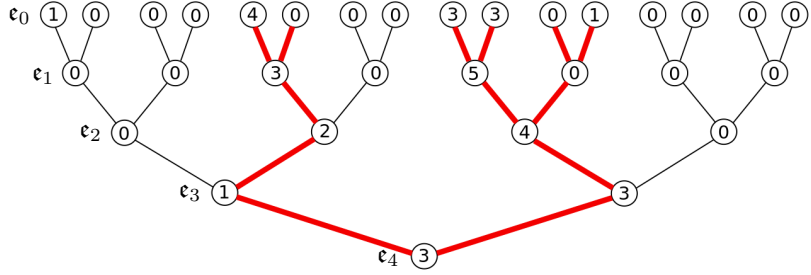

For $n\ge 0$, let $\mathfrak{e}_n$ denote the first lexicographic vertex in the $n$-th generation of $\T$. See Figure \ref{f:fig2}. We write
\begin{equation}
    N_n^\# := N^\# (\mathfrak{e}_n),
    \qquad
    N_n^{(i)} := N^{(i)} (\mathfrak{e}_n),
    \qquad
    Y_n := Y(\mathfrak{e}_n) \, .
    \label{N}
\end{equation}

\noindent [In the paper, only the law of each $Y_n$ {\it individually} is concerned; since $Y(\mathfrak{e}_n)$ has the same distribution as $Y_n$, the abuse of notation $Y_n := Y(\mathfrak{e}_n)$, which has the advantage of making formulae and discussions more compact, should not be source of any confusion.] See Figure \ref{f:fig2}.

\subsection{A bridge connecting two banks}
\label{subs:1st-key}

Assume $\P(X_0=k) \ge \P(Y_0 =k)$ for all integers $k\ge 1$, and $\P(Y_0=0)>0$. Let $(X_n, \, n\ge 0)$ and $(Y_n, \, n\ge 0)$ be systems coupled via the $XY$-coupling (see Theorem \ref{t:coupling}). Thanks to the notion of the number of open paths $N_n^{(0)}$ for $(Y_n)$ introduced in \eqref{N}, we are now able to give a lower bound for $\E(X_n)$ in terms of the number of open paths of $(Y_n)$.

\medskip

\begin{theorem}
\label{t:connexion}

 {\bf (The bridge inequality).}
 Assume $\P(X_0=k) \ge \P(Y_0 =k)$ for all integers $k\ge 1$, and $\P(Y_0=0)>0$. Let $(X_n, \, n\ge 0)$ and $(Y_n, \, n\ge 0)$ be systems coupled via the $XY$-coupling in Theorem \ref{t:coupling}. Let $\eta >0$. If $\E(X_0)<\infty$ and $\E(X_0-Y_0) \ge \eta \, \P(Y_0=0)$, then for all $r\ge 0$, all integers $n\ge 0$, $k\ge 0$ and $\ell \in [0, \, \frac{r \eta}{2}]$,
 \begin{equation}
     \E(X_{n+k+\ell})
     \ge
     \frac{m^{k+\ell} \, \eta}{2}\, \E\Big[ N_n^{(0)}\, {\bf 1}_{\{ N_n^{(0)} \ge r\} } \, {\bf 1}_{\{ Y_n =k\} } \Big] \, .
     \label{pf_lb_E(Xn)>_eq3}
 \end{equation}

\end{theorem}

\medskip

\noindent {\it Proof.} Let $A_n$ denote the set of vertices $x\in \T$ with $|x|=0$ that are the beginning of an open path leading to $\mathfrak{e}_n$. [So in the notation of \eqref{N}, the cardinality of $A_n$ is $N^\#_n$.] Let
$$
A_n^{(0)}
:=
\{ x\in A_n: Y(x) =0\} \, .
$$

\noindent The cardinality of $A_n^{(0)}$ is $N_n^{(0)}$, defined in \eqref{N}. Write
$$
\xi_n
:=
\sum_{x\in A_n^{(0)}} [X(x) - Y(x)] \, .
$$

\noindent [Recall that $X(x) \ge Y(x)$ a.s.\ and $\P(X(x) = Y(x) \, | \, Y(x) >0)=1$ for all $x\in \T$ with $|x| =0$.] The crucial observation is that for all $n\ge 1$,
$$
    X(\mathfrak{e}_n)
    \ge
    Y(\mathfrak{e}_n)
    +
    \xi_n ,
    \qquad \hbox{\rm a.s.}
$$

\noindent So for any integers $k\ge 0$ and $r\ge 0$, and any $\ell \ge 0$,
\begin{eqnarray*}
    (X(\mathfrak{e}_n)-k-\ell)^+
 &\ge& (Y(\mathfrak{e}_n) + \xi_n -k-\ell)^+ \, {\bf 1}_{\{ N_n^{(0)} \ge r\} } \, {\bf 1}_{\{ Y(\mathfrak{e}_n) =k\} }
    \\
 &=& (\xi_n -\ell)^+ \, {\bf 1}_{\{ N_n^{(0)} \ge r\} } \, {\bf 1}_{\{ Y(\mathfrak{e}_n) =k\} } \, .
\end{eqnarray*}

Let
$$
\mathscr{Y}
:=
\sigma( Y(x), \, x\in \T)
=
\sigma( Y(x), \, x\in \T, \, |x|=0) \, ,
$$

\noindent which is the sigma-field generated by the hierarchical system $(Y(x), \, x\in \T)$. Since $\{ N_n^{(0)} \ge r\} \in \mathscr{Y}$ and $\{ Y(\mathfrak{e}_n) =k\}\in \mathscr{Y}$, we obtain:
\begin{equation}
    \E ( (X(\mathfrak{e}_n)-k-\ell)^+ \, | \, \mathscr{Y} )
    \ge {\bf 1}_{\{ N_n^{(0)} \ge r\} } \, {\bf 1}_{\{ Y(\mathfrak{e}_n) =k\} } \, \E ( (\xi_n-\ell)^+ \, | \, \mathscr{Y} ) \, .
    \label{X>Y+eta}
\end{equation}

\noindent Conditionally on $x\in A_n^{(0)}$, the random variable $X(x) - Y(x)$ has the same law as $X_0-Y_0$ conditionally on $\{Y_0=0\}$; hence
$$
\E (\xi_n\, | \, \mathscr{Y})
=
N_n^{(0)} \, \E(X_0 - Y_0 \, | \, Y_0=0)
=
\frac{N_n^{(0)} \, \E(X_0 - Y_0)}{\P(Y_0=0)} \, ,
$$

\noindent the last equality being a consequence of \eqref{E(Z|Y)}. By assumption, this yields $\E (\xi_n\, | \, \mathscr{Y}) \ge \eta \, N_n^{(0)}$. Consequently, for any $\ell \ge 0$,
$$
    \E ( (\xi_n-\ell)^+ \, | \, \mathscr{Y} )
    \ge
    \E ( \xi_n -\ell \, | \, \mathscr{Y} )
    \ge
    \eta \, N_n^{(0)} - \ell \, .
$$

\noindent Combined with \eqref{X>Y+eta} and taking expectation, it follows that for any integers $k\ge 0$ and $r\ge 0$, and any $\ell \ge 0$,
$$
\E [(X(\mathfrak{e}_n)-k-\ell)^+]
\ge
\E \Big[ (\eta \, N_n^{(0)}- \ell) \, {\bf 1}_{\{ N_n^{(0)} \ge r\} } \, {\bf 1}_{\{ Y(\mathfrak{e}_n) =k\} } \Big] \, .
$$

\noindent On the event $\{ N_n^{(0)} \ge r\}$, we have $\eta \, N_n^{(0)} - \ell \ge \frac{\eta\, N_n^{(0)}}{2}$ if $\ell \le \frac{r \eta}{2}$. Hence, for any integers $k\ge 0$, $r\ge 0$ and $0\le \ell \le \frac{r \eta}{2}$,
$$
\E [(X(\mathfrak{e}_n)-k-\ell)^+]
\ge
\E \Big[ \frac{\eta}{2} \, N_n^{(0)}\, {\bf 1}_{\{ N_n^{(0)} \ge r\} } \, {\bf 1}_{\{ Y(\mathfrak{e}_n) =k\} } \Big] \, .
$$

\noindent If $\ell\ge 0$ is taken to be an integer, then obviously
$$
\E (X(\mathfrak{e}_{n+ k+\ell}))
\ge
m^{k+\ell} \, \E [(X(\mathfrak{e}_n)-k-\ell)^+] \, .
$$

\noindent Consequently, for any integers $k\ge 0$, $r\ge 0$ and $\ell \in [0, \, \frac{r \eta}{2}]$,
$$
\E (X(\mathfrak{e}_{n+ k+\ell}))
\ge
\frac{m^{k+\ell} \, \eta}{2}\, \E \Big[ N_n^{(0)}\, {\bf 1}_{\{ N_n^{(0)} \ge r\} } \, {\bf 1}_{\{ Y(\mathfrak{e}_n) =k\} } \Big] \, .
$$

\noindent Recall that $X(\mathfrak{e}_{n+ k+\ell})$ has the distribution of $X_{n+ k+\ell}$, and that $Y(\mathfrak{e}_n) =: Y_n$ by notation (see \eqref{N}). This yields Theorem \ref{t:connexion}.\qed

\subsection{The first crossing}
\label{subs:the_first_crossing}

Theorem \ref{t:connexion} yields the following estimate for the number of open paths in a critical system. This is our first crossing of the bridge, studying the critical system by means of a supercritical system. In Section \ref{s:2nd_crossing}, we are going to make a second crossing of the same bridge, but in the opposite direction, studying the supercritical system by means of a critical system.

\medskip

As in \eqref{Lambda}, we keep using the following notation:
$$
\Lambda = \Lambda(Y_0) := \E(Y_0^3\, m^{Y_0}) \in (0, \, \infty] \, ,
$$

\noindent associated with the critical system $(Y_n)$, i.e., satisfying $(m-1)\E(Y_0 \, m^{Y_0}) = \E(m^{Y_0}) <\infty$. Note that $\P(Y_0=0)>0$.\footnote{Recall that we always assume $\P(Y_0\ge 2)>0$ to avoid triviality.} Recall from \eqref{Lambda>} that for the critical system $(Y_n)$, we have
$$
\Lambda \in [\frac{1}{m-1} , \, \infty]\, .
$$

Let $\beta\in [0, \, 2]$ and $\chi \in (0, \, 1]$. We assume the system $(Y_n, \, n\ge 0)$ is $\beta$-regular with coefficient $\chi$ in the sense of \eqref{alpha}. Recall from Remark \ref{r:alpha=2;Lambda<infty} that $\Lambda<\infty$ in case $\beta=2$. Write
\begin{equation}
    \lambda_n
    =
    \lambda_n(Y_0, \, \beta)
    :=
    \min\{ \Lambda, \, n^\beta \} \, .
    \label{lambda_n}
\end{equation}

\medskip

\begin{corollary}
\label{c:E(N)>}

 Let $\beta\in [0, \, 2]$ and $\chi \in (0, \, 1]$.
 Assume $(m-1)\E(Y_0 \, m^{Y_0}) = \E(m^{Y_0}) <\infty$ and the system $(Y_n, \, n\ge 0)$ is $\beta$-regular with coefficient $\chi$ in the sense of \eqref{alpha}.
 Let $N_n^{(0)}$ be as in \eqref{N}.
 There exist constants $c_{11}>0$ and $c_{12}>0$, depending only on $m$ and $\beta$, such that for all $n\ge 2$,
 $$
 \E\Big[ m^{Y_n} N_n^{(0)}\, {\bf 1}_{\{ N_n^{(0)} \ge \frac{c_{11}}{\chi} \frac{n^2 \log n}{\lambda_n} \} } \, {\bf 1}_{\{ Y_n \le  n\} } \Big]
 \le
 \frac{c_{12}}{ \chi n^2} ,
 $$

 \noindent where $\lambda_n$ is defined in \eqref{lambda_n}.

\end{corollary}

\medskip

\noindent {\it Proof.} Fix $n\ge 2$. Let $c_2\in (0, \, 1]$ and $c_3 \in (0, \, 1]$ be the constants in Proposition \ref{p:ub}. Let $n_1=10n$ and $\eta=\frac{c_3c_2 \frac{\chi}{m+1}}{(m-1)^2} \, \frac{\lambda_{10n}}{(10n)^2}$. Then $(m-1)^2\eta \in (0, c_3 \chi) \subset (0, c_3)$.

Let $(X_n^{(\eta)}, \, n\ge 0)$ be a system with $\P(X_0^{(\eta)}=k)=\P(Y_0=k)$ for $k\ge 2$, and
\begin{eqnarray*}
    \P(X_0^{(\eta)}=0)
 &=& \max\{\P(Y_0=0)-\eta, \, 0\},
    \\
    \P(X_0^{(\eta)}=1)
 &=& \P(Y_0=1) + \min\{\P(Y_0=0), \, \eta\}.
\end{eqnarray*}

\noindent Then $(m-1)\E(X_0^{(\eta)}m^{X_0^{(\eta)}}) -\E(m^{X_0^{(\eta)}}) = (m-1)^2 \min\{\P(Y_0=0), \, \eta\}\le (m-1)^2\eta$. Moreover, $\E((X_0^{(\eta)})^3 m^{X_0^{(\eta)}}) = \E(Y_0^3 m^{Y_0}) + m \min\{\P(Y_0=0), \, \eta\} =\Lambda+ m \min\{\P(Y_0=0), \, \eta\}$. Since $\Lambda \ge \frac{1}{m-1}$, whereas $m \min\{\P(Y_0=0), \, \eta\}\le m \eta \le \frac{m}{(m-1)^2} \le \frac{m}{m-1} \le m \Lambda$, we get $\E((X_0^{(\eta)})^3 m^{X_0^{(\eta)}}) \le (m+1)\Lambda$. By assumption, $(Y_n, \, n\ge 0)$ is $\beta$-regular with coefficient $\chi$ in the sense of \eqref{alpha}, so the supercritical system $(X_n^{(\eta)}, \, n\ge 0)$ is also $\beta$-regular with coefficient $\frac{\chi}{m+1}$.

By definition of $n_1$ and $\eta$, we have $n_1 \le \min\{ (\frac{c_2 \, \frac{\chi}{m+1} \Lambda}{(m-1)^2\eta})^{1/2}, \, (\frac{c_2 \, \frac{\chi}{m+1}}{(m-1)^2\eta})^{1/(2-\beta)}\}$ for $\beta\in [0, \, 2)$ and $n_1 \le (\frac{c_2 \, \frac{\chi}{m+1}\Lambda}{(m-1)^2\eta})^{1/2}$ for $\beta=2$. So we are entitled to apply Proposition \ref{p:ub} to get
$$
\max_{0\le i\le n_1} \E(X_i^{(\eta)})\le 3.
$$

On the other hand, $\P(X_0^{(\eta)} =k) \ge \P(Y_0 =k)$ for all $k\ge 1$, and $\P(Y_0=0)>0$ as observed previously; since $\frac{\E(X_0^{(\eta)}-Y_0)}{\P(Y_0=0)} = \frac{\min\{\P(Y_0=0), \, \eta\}}{\P(Y_0=0)}\ge \eta$, we are entitled to apply Theorem \ref{t:connexion} to obtain, for all $r\ge 0$, all integers $k\ge 0$ and $\ell \in [0, \, \frac{r \eta}{2}]$,
$$
\E \Big[ N_n^{(0)}\, {\bf 1}_{\{ N_n^{(0)} \ge r\} } \, {\bf 1}_{\{ Y_n =k\} } \Big]
\le
\frac{2\E (X_{n+k+\ell}^{(\eta)}) }{m^{k+\ell} \, \eta} \, .
$$

\noindent Hence for all $r\ge 0$, and all integers $\ell \in [0, \, \frac{r \eta}{2}]$ with $2n+\ell\le n_1$, we have
$$
\E \Big[ m^{Y_n} N_n^{(0)}\, {\bf 1}_{\{ N_n^{(0)} \ge r\} } \, {\bf 1}_{\{ Y_n \le n\} } \Big]
\le
\frac{2(n+1)}{m^\ell \, \eta}  \max_{0\le i\le n_1} \E(X_i^{(\eta)}) \, \le \frac{6(n+1)}{m^\ell \, \eta}.
$$

\noindent We take $\ell := \lfloor  \frac{5\log n}{\log m}\rfloor$, $r:= \frac{c_{11}}{\chi} \frac{n^2 \log n}{\lambda_n}$ and $c_{11}:=\frac{10^3(m-1)^2 (m+1)}{c_3c_2 \log m}$. Then $2n+\ell\le 10n$ (because $\frac{5}{\log m} \le \frac{5}{\log 2} < 8$) and $\ell\le \frac{r\eta}{2}$. The corollary follows readily.\qed

\section{Preparation for the second crossing}
\label{s:open_paths}

Throughout the section, let $(Y_n, \, n\ge 0)$ denote a critical system, i.e., satisfying $(m-1)\E( Y_0 \, m^{Y_0}) = \E(m^{Y_0}) <\infty$.

The goal of this section is to establish a recursive formula (Proposition \ref{p:recursion_number_open_paths} below) for the weighted expected number of open paths.

The section is split into two parts. The first part gives a recursive formula for the number of open paths in the critical regime. The second part collects some known moment estimates for the system in the critical regime.

\subsection{A recursive formula}

Let $(Y_n, \, n\ge 0)$ be a system satisfying $(m-1)\E( Y_0 \, m^{Y_0}) = \E(m^{Y_0}) <\infty$. The following proposition gives a recursive formula for the weighted expected number of open paths, where $N_n$ can be $N_n^{(i)}$ (for any integer $i\ge 0$), or $N_n^\#$, as defined in \eqref{N}. Note that by definition, $\E[m^{Y_0}(1+Y_0)N_0^{(i)}]= m^i (1+i) \P(Y_0=i)$ and $\E[m^{Y_0}(1+Y_0)N_0^\#]=\E[m^{Y_0}(1+Y_0)]<\infty$.

In our paper, only the formula for $N_n=N_n^{(0)}$ is of interest.

\medskip

\begin{proposition}
 \label{p:recursion_number_open_paths}

 {\bf (Recursive formula for number of open paths).}
 Assume $(m-1)\E( Y_0 \, m^{Y_0}) = \E(m^{Y_0}) <\infty$. Fix an integer $i\ge 0$. Let $N_n$ denote either $N_n^\#$ or $N_n^{(i)}$. Then
 $$
 \E[m^{Y_n}(1+Y_n)N_n]
 =
 \E[m^{Y_0}(1+Y_0)N_0] \prod_{k=0}^{n-1} [\E(m^{Y_k})]^{m-1}\, ,
 \qquad
 n\ge 1 \, .
 $$

\end{proposition}

\medskip

\noindent {\it Proof.} For any vertex $x\in \T$ with $|x|=n \ge 1$, let $x^{(1)}$, $\ldots$, $x^{(m)}$ be as before its parents in generation $n-1$. By definition,
\begin{eqnarray*}
    Y(x)
 &=& [Y(x^{(1)}) + \cdots + Y(x^{(m)})-1]^+,
    \\
    N(x)
 &=& [N(x^{(1)}) + \cdots + N(x^{(m)})]\, {\bf 1}_{\{ Y(x^{(1)}) + \cdots + Y(x^{(m)}) \ge 1\} }  \, ,
\end{eqnarray*}

\noindent where $N(y)$ can be either $N^\#(y)$ or $N^{(i)}(y)$. Write $\Sigma_Y = \Sigma_Y(x) := Y(x^{(1)}) + \cdots + Y(x^{(m)})$ and $\Sigma_N = \Sigma_N(x) := N(x^{(1)}) + \cdots + N(x^{(m)})$. Then $(Y(x), \, N(x)) = (\Sigma_Y -1, \, \Sigma_N) \, {\bf 1}_{ \{ \Sigma_Y \ge 1\} } + (0, \, 0) \, {\bf 1}_{ \{ \Sigma_Y =0\} }$. Consequently, for $s>0$ and $t\ge 0$ (with $0^0 := 1$),
\begin{eqnarray}
    \E(s^{Y(x)} \, t^{N(x)})
 &=& \E(s^{\Sigma_Y -1} \, t^{\Sigma_N}\, {\bf 1}_{\{ \Sigma_Y \ge 1 \} })
    +
    \P(\Sigma_Y =0)
    \nonumber
    \\
 &=& \E(s^{\Sigma_Y -1} \, t^{\Sigma_N})
    -
    s^{-1} \, \E(t^{\Sigma_N}\, {\bf 1}_{\{ \Sigma_Y =0 \} })
    +
    \P(\Sigma_Y =0) .
    \label{pf:recursive_formula1}
\end{eqnarray}

\noindent For further use, we observe that the same argument yields $\P(N(x) \ge 1, \, Y(x) =0) = \P(\Sigma_N \ge 1, \, \Sigma_Y =1)$, which is bounded by $\P( \Sigma_Y \ge 1)$. Hence
\begin{equation}
    \P(N(x) \ge 1, \, Y(x) =0)
    \le
    m \, \P(Y(x^{(1)}) \ge 1) \, .
    \label{pf:recursive_formula2}
\end{equation}

Define
$$
G_n(s, \, t)
:=
\E(s^{Y_n} \, t^{N_n}) ,
\qquad
s\ge 0, \; t \ge 0,
$$

\noindent which is the joint moment generating function for the pair $(Y_n, \, N_n)$. Since $(Y(x^{(i)}), \, N(x^{(i)}))$, for $1\le i\le m$, are i.i.d.\ having the distribution of $(Y_{n-1}, \, N_{n-1})$, we have $\E(s^{\Sigma_Y -1} \, t^{\Sigma_N}) = s^{-1} \, G_{n-1} (s, \, t)^m$, $\E(t^{\Sigma_N}\, {\bf 1}_{\{ \Sigma_Y =0 \} }) = [\E(t^{N_{n-1}}\, {\bf 1}_{\{ Y_{n-1} =0 \} })]^m = G_{n-1}(0, \, t)^m$, and $\P(\Sigma_Y =0) = [\P( Y_{n-1} =0)]^m = G_{n-1}(0, \, 1)^m$. So \eqref{pf:recursive_formula1} reads: for $n\ge 1$, $s>0$ and $t\ge 0$,
\begin{equation}
    G_n (s, \, t)
    =
    \frac{G_{n-1}(s, \, t)^m}{s}
    -
    \frac{G_{n-1}(0, \, t)^m}{s}
    +
    G_{n-1}(0, \, 1)^m \, .
    \label{pf:recursive_formula3}
\end{equation}

For further use, we also rewrite \eqref{pf:recursive_formula2} as follows:
\begin{equation}
    \P(N_n \ge 1, \, Y_n =0)
    \le
    m \, \P(Y_{n-1}\ge 1) \, .
    \label{pf:recursive_formula4}
\end{equation}

Let
\begin{equation}
    I_n(s)
    :=
    \E(s^{Y_n} \, N_n) \, ,
    \label{In}
\end{equation}

\noindent which is the partial derivative of $G_n(s, \, t)$ with respect to $t$ at $1$. Then
\begin{equation}
    I_{n+1}(s)
    =
    \frac{m}{s} \, I_n(s) G_n(s)^{m-1}
    -
    \frac{m}{s} \, I_n(0) G_n(0)^{m-1} ,
    \label{e:In_recursion}
\end{equation}

\noindent where
\begin{equation}
    G_n(s)
    :=
    \E(s^{Y_n}) .
    \label{Gn}
\end{equation}

\noindent [So $G_n(s) = G_n(s, \, 1)$.] Multiplying by $s$ on both sides of \eqref{e:In_recursion}, and differentiating with respect to $s$, we get
\begin{equation}
    I_{n+1}(s) + s \, I_{n+1}'(s)
    =
    m(m-1) G_n(s)^{m-2} G_n'(s) I_n(s)
    +
    m G_n(s)^{m-1} I_n'(s).
    \label{e:NXgene}
\end{equation}

\noindent In the critical regime, we have $m(m-1)G_n'(m)=G_n(m)$ (see \eqref{crit}), which yields, for $n\ge 0$,
\begin{equation}
    I_{n+1}(m) + mI_{n+1}'(m)
    =
    G_n(m)^{m-1}(I_n(m) + mI_n'(m)) \, .
    \label{In_recursion}
\end{equation}

\noindent Iterating the identity, and noting that $\E[m^{Y_n}(1+Y_n)N_n] = I_n(m)+mI_n'(m)$, we get
$$
\E[m^{Y_n}(1+Y_n)N_n]
=
\E[m^{Y_0}(1+Y_0)N_0] \prod_{k=0}^{n-1} G_k(m)^{m-1} .
$$

\noindent Proposition \ref{p:recursion_number_open_paths} is proved.\qed

\subsection{Moment estimates}

The assumption $(m-1)\E( Y_0 \, m^{Y_0}) = \E(m^{Y_0}) <\infty$ implies (Collet et al.~\cite{collet-eckmann-glaser-martin}) that for all $n\ge 0$,
\begin{equation}
    (m-1)\E( Y_n \, m^{Y_n}) = \E(m^{Y_n}) < \infty \, .
    \label{crit}
\end{equation}

\noindent By \eqref{Gn<}, with $c_{13} := (1+\frac{1}{m-1})m^{1/(m-1)}$,
\begin{equation}
    \E[(1+Y_n)m^{Y_n}] \le c_{13},
    \qquad
    n\ge 0  \, .
    \label{moment_Y_under_Q}
\end{equation}

We keep using the notation introduced previously:
\begin{eqnarray*}
    \Lambda = \Lambda(Y_0)
 &:=& \E(Y_0^3 m^{Y_0}) \in [ \frac{1}{m-1}, \, \infty],
    \\
    \lambda_n = \lambda_n(Y_0)
 &:=& \min\{ \Lambda, \, n^\beta \} \, .
\end{eqnarray*}

\noindent The next is a collection of known useful properties concerning the moments of $Y_n$.

\medskip

\begin{fact}
 \label{f:prod_Gn}

 {\bf (\cite{bmxyz_questions})}
 Let $\beta\in [0, \, 2]$ and $\chi \in (0, \, 1]$.
 Assume $(m-1)\E( Y_0 \, m^{Y_0}) = \E(m^{Y_0})<\infty$. Assume the system $(Y_n, \, n\ge 0)$ is $\beta$-regular with coefficient $\chi$ in the sense of \eqref{alpha}. Then there exist constants $c_{14}>0$, $c_{15} >0$, $c_{16}>0$ and $c_{17}>0$, depending only on $m$ and $\beta$, such that for all $n\ge 1$,
\begin{eqnarray}
    \frac{c_{14}}{\Lambda} \, n^2
 &\le& \prod_{k=0}^{n-1} [\E (m^{Y_k})]^{m-1}
    \le
   \frac{c_{15}}{\chi}\, \frac{n^2}{\lambda_n} ,
    \label{prod_Gn}
    \\
    \E(Y_n^2 \, m^{Y_n})
 &\le& c_{16} \, (\frac{\Lambda}{\chi}\, \frac{n^2}{\lambda_n})^{1/2} \, ,
    \label{Gn_2nd_derivative}
    \\
    \E(Y_n^3 \, m^{Y_n})
 &\le& c_{17} \, \frac{\Lambda}{\chi}\, \frac{n^2}{\lambda_n} \, .
    \label{Gn_3rd_derivative}
\end{eqnarray}

\end{fact}
The first inequality in \eqref{prod_Gn} follows from Proposition 2 in \cite{bmxyz_questions}, whereas the second inequality follows from the proof of Proposition 1 in \cite{bmxyz_questions} but by replacing the last displayed inequality in the proof (saying that $\Delta_0(m-\frac{m}{n+2}) \ge \frac{\P(X_0=0)}{2(n+2)^2}$ in the notation of \cite{bmxyz_questions}) by the assumption of the $\beta$-regularity in the sense of \eqref{alpha}. Inequality \eqref{Gn_3rd_derivative} follows from (19) in \cite{bmxyz_questions} and the second inequality in \eqref{prod_Gn}. Finally, \eqref{Gn_2nd_derivative} is a consequence of \eqref{moment_Y_under_Q}, \eqref{Gn_3rd_derivative} and the Cauchy--Schwarz inequality. All the inequalities are valid even in case $\Lambda=\infty$.

\medskip

\begin{remark}
 \label{l:joint_moment_Yn_Nn}

 Consider a system $(Y_n, \, n\ge 0)$ satisfying $(m-1)\E( Y_0 \, m^{Y_0}) = \E(m^{Y_0})<\infty$. Assume $(Y_n, \, n\ge 0)$ is $\beta$-regular (for some $\beta\in [0, \, 2]$) with coefficient $\chi \in (0, \, 1]$ in the sense of \eqref{alpha}.

Combining Proposition \ref{p:recursion_number_open_paths} with \eqref{prod_Gn} yields that there exist constants $c_{18}>0$ and $c_{19}>0$, depending only on $m$, such that
 \begin{equation}
     c_{18} \, \frac{n^2\, \P(Y_0=0)}{\Lambda}
     \le
     \E[m^{Y_n}(1+Y_n)N_n^{(0)}]
     \le
     \frac{c_{19}}{\chi}\, \frac{n^2\, \P(Y_0=0)}{\lambda_n}\, ,
     \label{joint_moment_Yn_Nn}
 \end{equation}

 \noindent for all $n\ge 1$. This is valid without assuming $\Lambda$ to be finite.

 About the factor $\P(Y_0=0)$ on both sides of \eqref{joint_moment_Yn_Nn}: the assumption $(m-1)\E( Y_0 \, m^{Y_0}) = \E(m^{Y_0})<\infty$ implies $\P(Y_0=0)>0$; in case $m\ge 3$, we have, moreover, $\P(Y_0=0)\ge \frac{m-2}{m-1}$ (Lemma 3 in \cite{bmxyz_questions}).\qed
\end{remark}

\section{The second crossing and a general lower bound}
\label{s:2nd_crossing}

In Section \ref{subs:the_first_crossing}, we crossed the bridge (Theorem \ref{t:connexion}) for the first time to obtain an upper bound for the number of open paths in a critical system.
Now, when the recursive formula in Proposition \ref{p:recursion_number_open_paths} is established, we are going to cross again the bridge, but in the opposite direction, to obtain a lower bound for the expected value in a supercritical system. The second crossing will lead to a general lower bound for the free energy of certain systems that are not necessarily integer-valued.

Let $\beta\in [0, \, 2]$ and $\chi \in (0, \, 1]$. Let $(Y_n, \, n\ge 0)$ be a system satisfying $(m-1)\E(Y_0 \, m^{Y_0}) = \E(m^{Y_0}) <\infty$ and is $\beta$-regular with coefficient $\chi$ in the sense of \eqref{alpha}. Let as before,
\begin{eqnarray*}
    \Lambda = \Lambda(Y_0)
 &:=& \E(Y_0^3 m^{Y_0}) \in [ \frac{1}{m-1}, \, \infty],
    \\
    \lambda_n = \lambda_n(Y_0)
 &:=& \min\{ \Lambda, \, n^\beta \} \, .
\end{eqnarray*}

\noindent By Remark \ref{r:alpha=2;Lambda<infty}, we have $\Lambda<\infty$ in case $\beta=2$.

We start with a couple of lemmas.

\medskip

\begin{lemma}
\label{l:E(N)<}

 Let $\beta\in [0, \, 2]$ and $\chi \in (0, \, 1]$.
 Assume $(m-1)\E(Y_0 \, m^{Y_0}) = \E(m^{Y_0}) <\infty$ and the system $(Y_n, \, n\ge 0)$ is $\beta$-regular with coefficient $\chi$ in the sense of \eqref{alpha}.
 Let $N_n^{(0)}$ be as in \eqref{N}.
 There exists a constant $c_{20}>0$, depending only on $m$ and $\beta$, such that for all integer $n \ge 2$,
 \begin{equation}
     \E (m^{Y_n} N_n^{(0)} )
     \le \frac{c_{20}}{\chi} \, \frac{n^2}{\lambda_n} \Big( \frac1n + (\log n) \, [\E(m^{Y_n} \, {\bf 1}_{\{ Y_n\ge 1\} })+\E(m^{Y_{n-1}} \, {\bf 1}_{\{ Y_{n-1}\ge 1\} })]\Big).
     \label{E(N)<}
 \end{equation}

\end{lemma}

\medskip

\noindent {\it Proof.} Let $c_{11}>0$ be the constant in Corollary \ref{c:E(N)>}.  Write
$$
\E(m^{Y_n} N_n^{(0)})
\le
J_1(n) + J_2(n) + J_3 (n) \, ,
$$
\noindent where
\begin{eqnarray*}
    J_1(n)
 &:=& \E [ m^{Y_n} N_n^{(0)} \, {\bf 1}_{\{ N_n^{(0)} > \frac{c_{11}}{\chi} \frac{n^2 \log n}{\lambda_n}\} } \, {\bf 1}_{\{ Y_n < n\} } ] \, ,
    \\
    J_2(n)
 &:=& \E [ m^{Y_n} N_n^{(0)} \, {\bf 1}_{\{ Y_n \ge n\} } ] \, ,
    \\
    J_3(n)
 &:=& \E [ m^{Y_n} N_n^{(0)} \, {\bf 1}_{\{ N_n^{(0)} \le \frac{c_{11}}{\chi} \frac{n^2 \log n}{\lambda_n} \} } ]
    \le
 \frac{c_{11}}{\chi} \frac{n^2 \log n}{\lambda_n} \, \E [ m^{Y_n} \, {\bf 1}_{\{ N_n^{(0)} \ge 1\} } ] \, .
\end{eqnarray*}

By Corollary \ref{c:E(N)>}, $J_1(n) \le \frac{c_{12} }{\chi n^2}\le \frac{c_{12}}{\chi }n^{1-\beta}$.

For $J_2(n)$, things are simple, because $m^{Y_n} N_n^{(0)} \, {\bf 1}_{\{ Y_n \ge n\} } \le \frac{1}{n} \, Y_n m^{Y_n} N_n^{(0)}$, so by \eqref{joint_moment_Yn_Nn}, for $n\ge 1$,
$$
J_2(n)
\le
\frac{c_{19}}{\chi} \, \frac{n\, \P(Y_0=0)}{\lambda_n}
\le
\frac{c_{19}}{\chi} \, \frac{n}{\lambda_n}\, .
$$

For $J_3(n)$, we argue that $\E [ m^{Y_n} \, {\bf 1}_{\{ N_n^{(0)} \ge 1\} } ] \le \P( N_n^{(0)} \ge 1, \, Y_n=0) + \E [ m^{Y_n} \, {\bf 1}_{\{ Y_n \ge 1\} } ]$. Recall from \eqref{pf:recursive_formula4} that $\P( N_n^{(0)} \ge 1, \, Y_n=0) \le m \, \P(Y_{n-1} \ge 1)$, which is bounded by $\E [ m^{Y_{n-1}} \, {\bf 1}_{\{ Y_{n-1} \ge 1\} } ]$. Hence
$$
J_3(n)
\le
\frac{c_{11}}{\chi} \frac{n^2 \log n}{\lambda_n} \, [\E ( m^{Y_n} \, {\bf 1}_{\{ Y_n \ge 1\} } )
+
\E ( m^{Y_{n-1}} \, {\bf 1}_{\{ Y_{n-1} \ge 1\} } )].
$$

\noindent Assembling these pieces yields the desired result.\qed

\bigskip

For the next lemma, let us write $I_n(s) := \E(s^{Y_n} \, N_n^{(0)})$ as in \eqref{In}.

\medskip

\begin{lemma}
\label{sum_I_n}

 Let $\beta\in [0, \, 2]$ and $\chi \in (0, \, 1]$. Assume $(m-1)\E(Y_0 \, m^{Y_0}) = \E(m^{Y_0}) <\infty$ and the system $(Y_n, \, n\ge 0)$ is $\beta$-regular with coefficient $\chi$ in the sense of \eqref{alpha}.
 There exists a constant $c_{21}>0$, depending only on $m$ and $\beta$, such that for all integer $n \ge 2$,
 $$
 \sum_{\ell=0}^n I_\ell (m)
 \le
 \frac{c_{21}}{{\chi^2}} \frac{n^2 (\log n)^2}{\lambda_n} \, .
 $$
\end{lemma}

\medskip

\noindent {\it Proof.} By Lemma \ref{l:E(N)<},
\begin{eqnarray*}
    \sum_{\ell=0}^n I_\ell (m)
 &\le& \E(m^{Y_0}N_0^{(0)})+\E(m^{Y_1}N_1^{(0)})
    \\
 && +
    \sum_{\ell=2}^n \frac{c_{20}}{\chi} \Big( \frac{\ell}{\lambda_\ell} + \frac{\ell^2 \log \ell}{\lambda_\ell} \, [\E(m^{Y_\ell} \, {\bf 1}_{\{ Y_\ell\ge 1\} }) + \E(m^{Y_{\ell-1}} \, {\bf 1}_{\{ Y_{\ell-1}\ge 1\} })]\Big) .
\end{eqnarray*}

\noindent On the right-hand side, $N_0^{(0)}\le 1$ and $N_1^{(0)}<m$, so $\E(m^{Y_0}N_0^{(0)})+\E(m^{Y_1}N_1^{(0)}) \le \E(m^{Y_0}) + m\, \E(m^{Y_1}) \le (m+1) c_{13}$ (by \eqref{moment_Y_under_Q}). On the other hand, $\frac{\ell^2}{\lambda_\ell} \le \frac{n^2}{\lambda_n}$. Hence
\begin{eqnarray}
    \sum_{\ell=0}^n I_\ell (m)
 &\le& (m+1) c_{13} + \frac{c_{20}}{\chi} \, \frac{n^2}{\lambda_n} \sum_{\ell=2}^n \frac{1}{\ell}
    \nonumber
    \\
 && +
    \frac{c_{20}}{\chi} \, \frac{n^2 \log n}{\lambda_n} \sum_{\ell=2}^n [\E(m^{Y_\ell} \, {\bf 1}_{\{ Y_\ell\ge 1\} }) + \E(m^{Y_{\ell-1}} \, {\bf 1}_{\{ Y_{\ell-1}\ge 1\} })]
    \nonumber
    \\
 &\le& (m+1) c_{13} + \frac{c_{20}}{\chi} \, \frac{n^2 \log n}{\lambda_n} + \frac{2c_{20}}{\chi} \, \frac{n^2 \log n}{\lambda_n} \sum_{\ell=0}^n\E(m^{Y_\ell} \, {\bf 1}_{\{ Y_\ell\ge 1\} }).
    \label{esumIup}
\end{eqnarray}

\noindent On the other hand, noting that $G_\ell(m)\le m^{1/(m-1)}$ (see \eqref{Gn<}) and there exists $c_{22}=c_{22}(m)<\infty$ such that $\log x \ge c_{22}(x-1)$ for $x\in [1, \, m^{1/(m-1)}]$, we get
$$
\E(m^{Y_\ell} \, {\bf 1}_{\{ Y_\ell\ge 1\} })
\le
2 \E[(m^{Y_\ell}-1) \, {\bf 1}_{\{ Y_\ell\ge 1\} }]
=
2[\E(m^{Y_{\ell}})-1]\le \frac{2}{c_{22}} \log [G_\ell(m)].
$$

\noindent Hence
\begin{eqnarray*}
    \sum_{\ell=0}^n \E(m^{Y_\ell}\, {\bf 1}_{\{ Y_\ell\ge 1\} })
 &\le& \frac{2}{c_{22}} \sum_{\ell=0}^n \log [G_\ell(m)]
    \\
 &=& \frac{2}{(m-1)c_{22}}\log \Big( \prod_{\ell=0}^n G_\ell(m)^{m-1} \Big)
    \\
 &\le& \frac{2}{(m-1)c_{22}} \log \Big( \frac{c_{15}}{\chi}\, \frac{(n+1)^2}{\lambda_{n+1}}\Big),
\end{eqnarray*}

\noindent the last inequality being a consequence of \eqref{prod_Gn}. Note that $\log (\frac{c_{15}}{\chi}\, \frac{(n+1)^2}{\lambda_{n+1}}) \le c_{23} \, \frac{\log n}{\chi}$ for some constant $c_{23}>0$ depending only on $m$, $\beta$ and $c_{15}$ (thus only on $m$ and $\beta$ once $c_{15}$ is chosen in \eqref{prod_Gn}). Substituting this into \eqref{esumIup}, we get the desired result.\qed

\medskip

\begin{proposition}
\label{p:third_derivative_In}

 Let $\beta\in [0, \, 2]$ and $\chi \in (0, \, 1]$.
 Assume $(m-1)\E(Y_0 \, m^{Y_0}) = \E(m^{Y_0}) <\infty$ and the system $(Y_n, \, n\ge 0)$ is $\beta$-regular with coefficient $\chi$ in the sense of \eqref{alpha}.
 Let $N_n^{(0)}$ be as in \eqref{N}.
 Then there exists a constant $c_{24}>0$, depending only on $m$ and $\beta$, such that for all $n\ge 2$,
 $$
 \E[(1+Y_n)^3\, m^{Y_n} N_n^{(0)}]
 \le
 c_{24} \, \frac{\Lambda}{\chi^3} \, \frac{n^4 (\log n)^2}{\lambda_n^2} \, .
 $$
\end{proposition}

\medskip

\noindent {\it Proof.} Let $G_n(s) := \E(s^{Y_n})$ as in \eqref{Gn}. Let $I_n(s) := \E(s^{Y_n} \, N_n^{(0)})$ as in \eqref{In}. [In particular, $I_0(s) = \P(Y_0=0)$.] Recall the iteration formula \eqref{e:NXgene}:
$$
I_{n+1}(s) + s \, I_{n+1}'(s)
=
m(m-1) G_n(s)^{m-2} G_n'(s) I_n(s)
+
m G_n(s)^{m-1} I_n'(s).
$$

\noindent Differentiating twice on both sides with respect to $s$ and taking $s=m$, and using the identity $(m-1)mG'_n(m)=G_n(m)$ in the critical regime (see \eqref{crit}), we deduce that for some unimportant positive constants $c_{25} >0$, $\ldots$, $c_{29}>0$ and all $n\ge 1$,
\begin{equation}
    3I_{n+1}''(m) + m \, I_{n+1}'''(m)
    =
    G_n(m)^{m-1}\Upsilon_n
    +
    G_n(m)^{m-1} [3 I_n''(m) + m\, I_n'''(m)]\, ,
    \label{e:I_'''<G}
\end{equation}

\noindent where
\begin{eqnarray*}
    \Upsilon_n
 &:=& c_{25} \, I_n(m)
    +
    c_{26} \, \frac{G_n''(m)}{G_n(m)}\, I_n(m)
    +
    c_{27} \, \frac{G_n'''(m)}{G_n(m)}\, I_n(m)
    \\
    && \qquad +
    c_{28} \, I_n'(m)
    +
    c_{29} \, \frac{G_n''(m)}{G_n(m)}\, I_n'(m)
    \\
 &\le& c_{25} \, I_n(m)
    +
    c_{26} \, G_n''(m) \, I_n(m)
    +
    c_{27} \, G_n'''(m) \, I_n(m)
     \\
    && \qquad +
    c_{28} \, I_n'(m)
    +
    c_{29} \, G_n''(m) \, I_n'(m) \, .
\end{eqnarray*}

At this stage, it is convenient to introduce
$$
\mathcal{D}_n(m)
:=
m(m-1)G_n'''(m) + (4m-5)G_n''(m) + \frac{2(m-2)}{m^2(m-1)} \, G_n(m)\, .
$$

\noindent Obviously, $\mathcal{D}_n(m)\ge c_{30}(1 +G_n''(m)+G_n'''(m))$ for some constant $c_{30}>0$ depending only on $m$, so that with $c_{31} := c_{30}\max\{c_{25}, \, c_{26}, \, c_{27} \}$ and $c_{32} := c_{28}+c_{29}$,
\begin{eqnarray*}
    \Upsilon_n
 &\le& c_{31} I_n(m)\mathcal{D}_n(m)+c_{32} \max\{ 1, \, G_n''(m)\} \, I_n'(m)
    \\
 &\le& c_{31} I_n(m)\mathcal{D}_n(m)+c_{32} \max\{ 1, \, G_n''(m)\} \, [I_n(m)+ mI_n'(m)] \, .
\end{eqnarray*}

\noindent The interest of $\mathcal{D}_n(m)$ lies in its recursive relation $\mathcal{D}_{n+1}(m) = \mathcal{D}_n(m)\, G_n(m)^{m-1}$ (see \cite{bmxyz_questions}), which leads to
$$
\mathcal{D}_n(m)
=
\mathcal{D}_0(m)\prod_{i=0}^{n-1}G_i(m)^{m-1}
\le
c_{33}\Lambda \prod_{i=0}^{n-1}G_i(m)^{m-1} \, ,
$$

\noindent where $c_{33}>0$ is a constant depending only on $m$. On the other hand, by \eqref{In_recursion},
$$
I_n(m)+mI_n'(m)
=
(I_0(m)+mI_0'(m)) \prod_{i=0}^{n-1}G_i(m)^{m-1}\, .
$$

\noindent Since $I_0(s) = \P(Y_0=0)$ for all $s$, we have $I_0(m)+mI_0'(m) = \P(Y_0=0) \le 1$, which implies $I_n(m)+mI_n'(m) \le \prod_{i=0}^{n-1}G_i(m)^{m-1}$. Consequently,
$$
\Upsilon_n
\le
\Big( c_{34}\Lambda I_n(m)+c_{34} \max\{ 1, \, G_n''(m)\} \Big) \prod_{i=0}^{n-1}G_i(m)^{m-1} ,
$$

\noindent with $c_{34} := \max\{ c_{31}c_{33}, \, c_{32}\}$. So we can rewrite \eqref{e:I_'''<G} as
\begin{eqnarray*}
 &&3 I_{n+1}''(m) + m \, I_{n+1}'''(m)
    \\
 &\le& c_{34}(\Lambda \, I_n(m)+ \max\{ 1, \, G_n''(m)\} ) \prod_{i=0}^{n}G_i(m)^{m-1}
  \\
  && \qquad   +
    G_n(m)^{m-1}(3 I_n''(m) + m \, I_n'''(m)) .
\end{eqnarray*}

\noindent Iterating the inequality, and noting that $I_0''(m) = I_0'''(m) =0$ (because $I_0(s) = \P(Y_0=0)$ for all $s$), we obtain:
$$
3I_n''(m) + mI_n'''(m)
\le
c_{34}\Big( \prod_{k=0}^{n-1} G_k(m)^{m-1}\Big) \sum_{\ell=0}^{n-1}( \Lambda \, I_\ell(m)+ \max\{ 1, \, G_\ell''(m)\} ) \, .
$$

\noindent By \eqref{prod_Gn}, $\prod_{k=0}^{n-1} G_k(m)^{m-1} \le \frac{c_{15}}{\chi} \, \frac{n^2}{\lambda_n}$, whereas by Lemma \ref{sum_I_n}, $\sum_{\ell=0}^{n-1} I_\ell (m)\le \frac{c_{21}}{{\chi^2}} \frac{n^2 (\log n)^2}{\lambda_n}$; hence
$$
3I_n''(m) + mI_n'''(m)
\le
\frac{c_{34}c_{15}}{\chi}\, \frac{n^2}{\lambda_n} \Big( \Lambda \frac{c_{21}}{{\chi^2}} \frac{n^2 (\log n)^2}{\lambda_n} + \sum_{\ell=0}^{n-1} \max\{ 1, \, G_\ell''(m)\} \Big) \, .
$$

\noindent By \eqref{Gn_2nd_derivative}, $G_\ell''(m)\le \frac{c_{16}}{\chi^{1/2}} \, (\frac{\ell^2}{\lambda_\ell} \Lambda)^{1/2} = \frac{c_{16}}{\chi^{1/2}} \, \max \{ \ell, \, \Lambda^{1/2} \ell^{1-(\beta/2)}\} \le \frac{c_{16}}{\chi^{1/2}} \, (n+\Lambda^{1/2} n^{1-(\beta/2)})$; thus
\begin{eqnarray*}
    \sum_{\ell=0}^{n-1} \max\{ 1, \, G_\ell''(m)\}
 &\le& \sum_{\ell=0}^{n-1} \Big( 1 + \frac{c_{16}}{\chi^{1/2}} \, (n+\Lambda^{1/2} n^{1-(\beta/2)}) \Big)
    \\
 &=& n+ \frac{c_{16}}{\chi^{1/2}} \, (n^2+\Lambda^{1/2} n^{2-(\beta/2)})
    \\
 &\le& \frac{1+c_{16}}{\chi^{1/2}} \, (n^2+\Lambda^{1/2} n^{2-(\beta/2)})\, .
\end{eqnarray*}

\noindent Since $n^2+\Lambda^{1/2} n^{2-(\beta/2)} \le 2\max\{\Lambda^{1/2}n^{2-(\beta/2)}, \, n^2\} \le 2\max\{\Lambda n^{2-\beta}, \, n^2\} = 2\frac{n^2}{\lambda_n} \Lambda$, this yields
$$
\sum_{\ell=0}^{n-1} \max\{ 1, \, G_\ell''(m)\}
\le
\frac{2(1+c_{16})}{\chi^{1/2}} \, \frac{n^2}{\lambda_n} \Lambda \, .
$$

\noindent Consequently,
$$
3I_n''(m) + mI_n'''(m)
\le
\frac{c_{34}c_{15}}{\chi}\, \frac{n^2}{\lambda_n} \Big( \Lambda \frac{c_{21}}{{\chi^2}} \frac{n^2 (\log n)^2}{\lambda_n} + \frac{2(1+c_{16})}{\chi^{1/2}} \, \frac{n^2}{\lambda_n} \Lambda \Big) ,
$$

\noindent proving Proposition \ref{p:third_derivative_In} (recalling that $\chi \le 1$).\qed

\medskip

\begin{proposition}
\label{p:key_lb}

 Let $\beta\in [0, \, 2]$ and $\chi \in (0, \, 1]$. Assume $\Lambda := \E(Y_0^3 \, m^{Y_0}) <\infty$. Assume $(m-1)\E(Y_0 \, m^{Y_0}) = \E(m^{Y_0})$ and the system $(Y_n, \, n\ge 0)$ is $\beta$-regular with coefficient $\chi$ in the sense of \eqref{alpha}. Let $N_n^{(0)}$ be as in \eqref{N}.
 There exist positive constants $c_{35}$, $c_{36}$ and $c_{37}$ which depend only on $(m, \, \beta)$ such that, for all $n\ge 2$,
 \begin{eqnarray}
    && \E [(1+Y_n) m^{Y_n} N_n^{(0)} \, {\bf 1}_{\{ N_n^{(0)} > c_{35} \frac{n^2\, \P(Y_0=0)}{\Lambda}\} } \, {\bf 1}_{\{ Y_n < \frac{c_{36}}{\chi^{3/2}} \frac{n \, \Lambda \log n}{\lambda_n\, [\P(Y_0=0)]^{1/2}}\} } ]
    \\
    && \ge
     c_{37}\frac{n^2\, \P(Y_0=0)}{\Lambda}\, ,
     \label{key_lb_preparation} \nonumber
 \end{eqnarray}
 \noindent Consequently, there exist a constant $c_{38}>0$ depending only on $(m, \, \beta)$,
 and an integer $k$ with $0\le k<\frac{c_{36}}{\chi^{3/2}} \frac{n \, \Lambda \log n}{\lambda_n\, [\P(Y_0=0)]^{1/2}}$, such that
 \begin{equation}
     \E \Big[ N_n^{(0)} \, {\bf 1}_{\{ N_n^{(0)} > c_{35} \frac{n^2\, \P(Y_0=0)}{\Lambda}\} } \, {\bf 1}_{\{ Y_n=k\} } \Big]
     \ge
     c_{38}\, \frac{\chi^3}{m^k\Lambda^3}\, \frac{\lambda_n^2\, [\P(Y_0=0)]^2}{(\log n)^2} \, .
     \label{key_lb}
 \end{equation}

\end{proposition}

\medskip

\noindent {\it Proof.} Only \eqref{key_lb_preparation} needs to be proved. Let $c_{39}>0$ and let $K\ge 1$ be an integer; their values will be chosen later. By Proposition \ref{p:third_derivative_In},
$$
\E [(1+Y_n) m^{Y_n} N_n^{(0)} \, {\bf 1}_{\{ Y_n\ge K\} } ]
\le
\frac{\E [(1+Y_n)^3 m^{Y_n} N_n^{(0)}]}{K^2}
\le
\frac{c_{24}}{\chi^3K^2} \, \frac{n^4 \Lambda (\log n)^2}{\lambda_n^2}\, .
$$

\noindent On the other hand, $N_n^{(0)} \, {\bf 1}_{\{ N_n^{(0)} \le c_{39} \frac{n^2\, \P(Y_0=0)}{\Lambda} \} } \le c_{39} \frac{n^2\, \P(Y_0=0)}{\Lambda}$, so by \eqref{moment_Y_under_Q},
\begin{eqnarray*}
    \E \Big[(1+Y_n) m^{Y_n} N_n^{(0)} \, {\bf 1}_{\{ N_n^{(0)} \le c_{39} \frac{n^2\, \P(Y_0=0)}{\Lambda}\} } \Big]
 &\le& c_{39} \frac{n^2\, \P(Y_0=0)}{\Lambda}\, \E [(1+Y_n) m^{Y_n}]
    \\
 &\le& c_{13}  c_{39} \frac{n^2\, \P(Y_0=0)}{\Lambda}\, .
\end{eqnarray*}

\noindent Recall from \eqref{joint_moment_Yn_Nn} that $\E [(1+Y_n) m^{Y_n} N_n^{(0)} ] \ge \frac{c_{18} n^2\, \P(Y_0=0)}{\Lambda}$. Hence
\begin{eqnarray*}
 &&\E [(1+Y_n) m^{Y_n} N_n^{(0)} \, {\bf 1}_{\{ N_n^{(0)} > c_{39} \frac{n^2\, \P(Y_0=0)}{\Lambda}\} } \, {\bf 1}_{\{ Y_n < K\} } ]
    \\
 &\ge& \E [(1+Y_n) m^{Y_n} N_n^{(0)} ]
    -
    \E [(1+Y_n) m^{Y_n} N_n^{(0)} \, {\bf 1}_{\{ Y_n \ge K\} } ]
    \\
 && \qquad
    -
    \E [(1+Y_n) m^{Y_n} N_n^{(0)} \, {\bf 1}_{\{ N_n^{(0)} \le c_{39} \frac{n^2\, \P(Y_0=0)}{\Lambda}\} } ]
    \\
 &\ge& \frac{c_{18} n^2\, \P(Y_0=0)}{\Lambda}- \frac{c_{24}}{\chi^3K^2} \,\frac{n^4 \Lambda (\log n)^2}{\lambda_n^2} - c_{13} c_{39} \frac{n^2\, \P(Y_0=0)}{\Lambda} \, .
\end{eqnarray*}

\noindent Choosing $c_{39} := \frac{c_{18}}{3c_{13}}$ and $K:= \lfloor (\frac{3c_{24}}{c_{18}})^{1/2} \chi^{-3/2} \frac{n \, \Lambda \log n}{\lambda_n\, [\P(Y_0=0)]^{1/2}} \rfloor +1$ yields \eqref{key_lb_preparation}. Proposition \ref{p:key_lb} is proved.\qed

\bigskip

We are now ready for the second crossing of the bridge which is Theorem \ref{t:connexion}. It implies a general lower bound for the free energy of certain supercritical systems $(X_n, \, n\ge 0)$. We are also going to apply this general result to prove the lower bounds in Theorem \ref{t:main} (Section \ref{s:lb}) and in Theorems \ref{t:2<alpha<4} and \ref{t:alpha=2} (Section \ref{s:alpha}).

Assume $(m-1)\E(Y_0 \, m^{Y_0}) = \E(m^{Y_0}) <\infty$. In particular, $\P(Y_0=0)>0$. Assume $\P(X_0=k) \ge \P(Y_0 =k)$ for all integers $k\ge 1$. We can couple two systems $(X_n, \, n\ge 0)$ and $(Y_n, \, n\ge 0)$ via the $XY$-coupling in Theorem \ref{t:coupling}.

\medskip

\begin{theorem}
\label{t:2nd_crossing}

 Let $\beta\in [0, \, 2]$ and $\chi \in (0, \, 1]$. Assume $\Lambda = \Lambda(Y_0) := \E(Y_0^3 \, m^{Y_0})<\infty$. Assume $(m-1)\E(Y_0 \, m^{Y_0}) = \E(m^{Y_0})$ and the system $(Y_n, \, n\ge 0)$ is $\beta$-regular with coefficient $\chi$ in the sense of \eqref{alpha}. Assume $\P(X_0=k) \ge \P(Y_0 =k)$ for all integers $k\ge 1$. Let $(X_n, \, n\ge 0)$ and $(Y_n, \, n\ge 0)$ be systems coupled via the $XY$-coupling in Theorem \ref{t:coupling}. There exists a constant $c_{40}>0$, depending only on $m$ and $\beta$, such that if $\E(X_0-Y_0) \ge \eta$ for some $\eta\in (0, \, \frac{1}{3(m-1)}]$,\footnote{The choice of $\frac{1}{3(m-1)}$ is to ensure $R:= \frac{\Lambda}{\eta \, \chi\, \P (Y_0=0)} \ge 3$.} then
 \begin{equation}
     \max_{0\le j\le n_0 } \E (X_j) \ge 2\, ,
     \label{lb_eq1}
 \end{equation}
 where $n_0:=\lfloor \frac{c_{40}}{[\,\chi^3\, \P(Y_0=0)]^{1/2}} \max\{\Lambda (R \log R)^{(1-\beta)/2}, \, (R\log R)^{1/2}\} \log R\rfloor$ with $R:= \frac{\Lambda}{\eta \, \chi\, \P (Y_0=0)}$; in particular, the free energy $F_\infty^X$ of the system $(X_n, \, n\ge 0)$ satisfies
 \begin{equation}
     F_\infty^X
     \ge
     \exp \Big[ - \frac{c_{41}}{[\,\chi^3\, \P(Y_0=0)]^{1/2}} \max\{\Lambda (R\log R)^{(1-\beta)/2}, \, (R \log R)^{1/2}\} \log R \Big] \, ,
     \label{lb_eq2}
 \end{equation}
 where $c_{41}>0$ is a constant whose value depends only on $m$ and $\beta$.

\end{theorem}

\medskip

\noindent {\it Proof.} Only \eqref{lb_eq1} needs proving, because \eqref{lb_eq2} follows immediately from \eqref{lb_eq1} by means of \eqref{F_encadrement}.

Let $((X(x), \, Y(x)), \, x\in \T)$ be the pair of systems in the $XY$-coupling in Theorem \ref{t:coupling} satisfying $\E(X_0-Y_0) \ge \eta$. Recall from Theorem \ref{t:connexion} that for all $r\ge 0$, all integers $n\ge 3$, $k\ge 0$ and $\ell \in [0, \, \frac{r \eta}{2\, \P(Y_0=0)}]$,\footnote{Theorem \ref{t:connexion} only requires $n\ge 0$; we take $n\ge 3$ so that $\log n \ge \log 3 >1$.}
$$
\E(X_{n+k+\ell})
\ge
\frac{m^{k+\ell} \, \eta}{2\, \P(Y_0=0)}\, \E\Big[ N_n^{(0)}\, {\bf 1}_{\{ N_n^{(0)} \ge r\} } \, {\bf 1}_{\{ Y_n =k\} } \Big] \, .
$$

  By \eqref{key_lb} in Proposition \ref{p:key_lb}, there exists an integer $k$ with $0\le k<\frac{c_{36}}{\chi^{3/2}} \frac{n \, \Lambda \log n}{\lambda_n\, [\P(Y_0=0)]^{1/2}}$, such that
$$
\E \Big[ N_n^{(0)} \, {\bf 1}_{\{ N_n^{(0)} > c_{35} \frac{n^2\, \P(Y_0=0)}{\Lambda}\} } \, {\bf 1}_{\{ Y_n=k\} } \Big]
\ge
c_{38}\, \frac{\chi^3}{m^k\Lambda^3}\, \frac{\lambda_n^2\, [\P(Y_0=0)]^2}{(\log n)^2} \, ,
$$

\noindent which implies that for some $0\le k<\frac{c_{36}}{\chi^{3/2}} \frac{n \, \Lambda \log n}{\lambda_n\, [\P(Y_0=0)]^{1/2}}$ and all $\ell \in [0, \, \frac{c_{35}n^2 }{2\Lambda}\eta]$,
\begin{eqnarray*}
\E (X_{n+k+\ell})
&\ge&
\frac{m^{k+\ell} \, \eta}{2\, \P(Y_0=0)}\, c_{38}\, \frac{\chi^3}{m^k\Lambda^3}\, \frac{\lambda_n^2\, [\P(Y_0=0)]^2}{(\log n)^2}
\\
&=&
\frac{c_{38}\chi^3}{2}\, \frac{m^\ell \, \eta}{\Lambda^3}\, \frac{\lambda_n^2\, \P(Y_0=0)}{(\log n)^2} \, .
\end{eqnarray*}

\noindent Since $\lambda_n^2 = \min\{ \Lambda^2, \, n^{2\beta} \} \ge \frac{1}{(m-1)^2}$ (recalling that $\Lambda \ge \frac{1}{m-1}$ for critical systems), this yields, with $c_{42} := \frac{c_{38}}{2(m-1)^2}$,
\begin{eqnarray*}
    \E (X_{n+k+\ell})
 &\ge& c_{42}\, \chi^3 \, \frac{m^\ell \, \eta\, \P(Y_0=0)}{\Lambda^3 (\log n)^2}
    \\
 &\ge& \frac{c_{42} \, m^\ell}{(\log n)^2} \, \Big( \frac{\eta \, \chi\, \P (Y_0=0)}{\Lambda} \Big)^{\! 3}
    \\
 &=& \frac{c_{42} \, m^\ell}{R^3 (\log n)^2} \, .
\end{eqnarray*}

\noindent We take $\ell := \lfloor c_{43} (\log R + \log \log n) \rfloor$, where $c_{43}>0$ is a constant depending only on $m$ and $\beta$, such that $\frac{c_{42} \, m^\ell}{R^3 (\log n)^2} \ge 2$. On the other hand, the condition $\ell \le \frac{c_{35}n^2 }{2\Lambda}\eta$ is respected if we take $n:= \lfloor c_{44} (R\log R)^{1/2}\rfloor$ for some appropriate constant $c_{44}\ge 1$ whose value depends only on $m$ and $\beta$; moreover, $c_{44}\ge 1$ is taken to be sufficiently large so that $\log R \le n$. Accordingly,
$$
\E (X_{n+k+\ell}) \ge 2\, .
$$

\noindent By definition, $n\le \frac{n \Lambda}{\lambda_n} \le \chi^{-3/2} \frac{n \Lambda \log n}{\lambda_n}$; also, $\log R \le n$ and $\log \log n \le n$, thus $\log R + \log \log n \le 2n \le 2\chi^{-3/2} \frac{n \Lambda \log n}{\lambda_n}$, which yields $\ell := \lfloor c_{43} (\log R + \log \log n) \rfloor \le 2c_{43} \, \chi^{-3/2} \frac{n \Lambda \log n}{\lambda_n}$. And of course, $k< c_{36} \, \chi^{-3/2} \frac{n \Lambda \log n}{\lambda_n\, [\P(Y_0=0)]^{1/2}}$ by definition. As a result, with $c_{45} := 1+ 2c_{43} + c_{36}$, we have
$$
n+k+\ell
\le
\frac{c_{45}}{\chi^{3/2}} \frac{n \Lambda \log n}{\lambda_n\, [\P(Y_0=0)]^{1/2}}
=
\frac{c_{45}}{[\,\chi^3\, \P(Y_0=0)]^{1/2}} \max \{ \Lambda n^{1-\beta}, \, n\} \log n ,
$$

\noindent with $n:= \lfloor c_{44} (R \log R)^{1/2}\rfloor$. This yields \eqref{lb_eq1}.\qed

\section{Proof of Theorem \ref{t:main}: lower bound}
\label{s:lb}

Let $(X_n, \, n\ge 0)$ be a system satisfying $\E[(X_0^*)^3 \, m^{X_0^*}] <\infty$ and let $p=p_c+\varepsilon$. Let $(Y_n, \, n\ge 0)$ be the critical system coupled with $(X_n, \, n\ge 0)$ as in the $XY$-coupling in Theorem \ref{t:coupling}. Then $X_0 \ge Y_0$ a.s., $\E(X_0-Y_0) = (p-p_c) \, \E(X_0^*) = \varepsilon \, \E(X_0^*)$, $\Lambda = \Lambda(Y_0) := \E(Y_0^3 m^{Y_0}) = p_c \, \E[(X_0^*)^3 \, m^{X_0^*}] <\infty$ and $\P (Y_0=0) = 1-p_c$. By Lemma \ref{A_4P(X_0=0)}, $(Y_n, \, n\ge 0)$ is $0$-regular with coefficient $\chi := 1-p_c$ in the sense of \eqref{alpha}. It follows from Theorem \ref{t:2nd_crossing} that there exists a constant $c_{46}>0$ such that for all sufficiently small $\varepsilon>0$, the free energy $F_\infty (p_c+\varepsilon)$ of the system $(X_n, \, n\ge 0)$ satisfies
\begin{equation}
    F_\infty (p_c+\varepsilon)
    \ge
    \exp \Big\{ - c_{46} \, \frac{[\log ( \frac{1}{\varepsilon} )]^{3/2}}{\varepsilon^{1/2}}\, \Big\} \, ,
    \label{lb}
\end{equation}

\noindent which readily yields the lower bound in Theorem \ref{t:main}.\qed

\section{Proof of Theorems \ref{t:2<alpha<4} and \ref{t:alpha=2}}
\label{s:alpha}

We start with a simple comparison result, which is useful in the proof of the upper bound in Theorems \ref{t:2<alpha<4} and \ref{t:alpha=2}.

\medskip

\begin{lemma}
\label{l:comparison}

 Let $(U_n, \, n\ge 0)$ and $(V_n, \, n\ge 0)$ be two recursive systems. If $U_0 \ge V_0$ \hbox{\rm a.s.} and if $\E(U_0) <\infty$, then for all $n\ge 0$,
 $$
 0\le \E(U_n)-\E(V_n) \le m^n \, [\E(U_0)-\E(V_0)] \, .
 $$

\end{lemma}

\medskip

\noindent {\it Proof.} Since $U_0 \ge V_0$ a.s., we can couple the two systems so that $U_n \ge V_n$ a.s.\ for all $n\ge 0$; in particular, $\E(U_n)-\E(V_n) \ge 0$.

Let $(X_n, \, n\ge 0)$ be an arbitrary recursive system. Recall from the recurrence relation \eqref{Gn_iteration} that $\E (s^{X_{n+1}}) = \frac1s \, [\E (s^{X_n})]^m  + (1-\frac1s) \, [\P (X_n=0)]^m$. Differentiating with respect to $s$, and taking $s=1$, this yields $\E (X_{n+1}) = m \, \E(X_n) -1 + [\P (X_n=0)]^m$ for all $n\ge 0$. Iterating the identity yields that for $n\ge 0$,
$$
\E (X_n)
=
m^n \, \E(X_0)
-
\sum_{i=0}^{n-1} m^{n-i-1} \, \{ 1- [\P (X_i=0)]^m\} \, .
$$

\noindent Applying the identity to $X_n=U_n$ and to $X_n=V_n$, and taking the difference, we obtain:
$$
\E(U_n)-\E(V_n)
=
m^n \, [\E(U_0)-\E(V_0)]
+
\sum_{i=0}^{n-1} m^{n-i-1} \, \{ [\P (U_i=0)]^m - [\P (V_i=0)]^m\} \, .
$$

\noindent By the coupling, $U_i \ge V_i$ a.s.\ so $[\P (U_i=0)]^m - [\P (V_i=0)]^m \le 0$. The lemma follows.\qed

\bigskip

Theorems \ref{t:2<alpha<4} and \ref{t:alpha=2} are proved by means of a truncation argument. Let $X_0^*$ be a random variable taking values in $\{1, 2, \ldots\}$ and satisfying
\begin{equation}
    \P (X_0^* = k)
    \, \sim \,
    c_0 \, m^{-k} k^{-\alpha} \, ,
    \qquad
    k\to \infty \, ,
    \label{initial_law}
\end{equation}

\noindent for some $0<c_0<\infty$ and $2\le \alpha \le 4$ (in Theorem \ref{t:2<alpha<4}, we assume $2<\alpha\le 4$, whereas in Theorem \ref{t:alpha=2}, we assume $\alpha=2$). Let $\ell_0 \ge 1$ be the smallest integer $k\ge 1$ such that $\P (X_0^* = k)>0$. Let $M>\ell_0$ be an integer (which will ultimately tend to infinity). Let $\varrho = \varrho(M, \, \alpha)$ be a random variable taking values in $\{ 0, \, \ell_0\}$; we are going to give the distribution of $\varrho$ later. Consider the truncated random variable
\begin{equation}
    X_0^{(M)}
    :=
    \varrho \, {\bf 1}_{\{ X_0^* = \ell_0\} }
    +
    X_0^* \, {\bf 1}_{\{ \ell_0< X_0^* \le M\} } \, .
    \label{X0M}
\end{equation}

\noindent Let $Y_0^{(M)}$ be a random variable whose distribution is given by
$$
P_{Y_0^{(M)}}
=
p_M \, P_{X_0^{(M)}}
+
(1-p_M) \delta_0 \, ,
$$

\noindent where $p_M \in (0, \, 1)$ is such that
$$
\E_{p_M} (m^{Y_0^{(M)}})
=
(m-1) \, \E_{p_M} (Y_0^{(M)} \, m^{Y_0^{(M)}}) \, .
$$

\noindent [The subscript $p_M$ is to indicate the weight $p_M$ of $X_0^{(M)}$ in the distribution of $Y_0^{(M)}$.] The existence of $p_M \in (0, \, 1)$ is simple; actually the value of $p_M$ can be explicitly computed:
\begin{equation}
    p_M = \frac{1}{1+ \E \{ [(m-1) X_0^{(M)} -1] m^{X_0^{(M)}}\} } \in (0, \, 1) \, .
    \label{p_M}
\end{equation}

\noindent Hence
\begin{eqnarray}
    \frac{1}{p_M} -1
 &=& \E \{ [(m-1) X_0^{(M)} -1] m^{X_0^{(M)}}\}
    \nonumber
    \\
 &=& \E \{ [(m-1) X_0^* -1] m^{X_0^*} \, {\bf 1}_{\{ X_0^* \le M\} } \}
    -
    \P(X_0^*>M)
    -
    a_M \, ,
    \label{1/p_M}
\end{eqnarray}

\noindent where
\begin{equation}
    a_M
    :=
    \E \Big[ \{ [(m-1)\ell_0 -1]m^{\ell_0} - [(m-1)\varrho -1]m^\varrho \} \, {\bf 1}_{\{ X_0^* = \ell_0\} } \Big]  \, .
    \label{a_M}
\end{equation}

\noindent In case $2<\alpha\le 4$, we have $p_c(X_0^*)>0$, and by \eqref{p_c},
$$
\frac{1}{p_c(X_0^*)} -1
=
\E \{ [(m-1) X_0^* -1] m^{X_0^*} \} \, ,
$$

\noindent which, in view of \eqref{1/p_M}, yields that for $2<\alpha\le 4$,
\begin{equation}
    \frac{1}{p_c(X_0^*)} - \frac{1}{p_M}
    =
    \P(X_0^*>M)
    +
    a_M
    +
    \E \{ [(m-1) X_0^* -1] m^{X_0^*} \, {\bf 1}_{\{ X_0^* > M\} } \} \, .
    \label{1/pc-1/pM}
\end{equation}

\subsection{Proof of Theorems \ref{t:2<alpha<4} and \ref{t:alpha=2}: upper bound}

To prove the upper bound, the choice of $\varrho$ in \eqref{X0M} is simple: we choose $\varrho=X_0^*$ (so the truncated random variable $X_0^{(M)}$ in \eqref{X0M} is simply $X_0^* \, {\bf 1}_{\{ X_0^* \le M\} }$). As such, $a_M=0$ in \eqref{a_M}.

The next lemma gives the asymptotic behaviour of $p_M-p_c(X_0^*)$ when $M\to \infty$. Note that $p_c(X_0^*)=0$ in case $\alpha=2$.

\medskip

\begin{lemma}
\label{l:p_M-p_c}

 Assume $\P(X_0^* = k) \sim c_0 \, m^{-k} k^{-\alpha}$, $k\to \infty$, for some $0<c_0<\infty$ and $2\le \alpha \le 4$. When $M\to \infty$,
 \begin{equation}
     p_M-p_c(X_0^*)
     \, \sim \,
     \begin{cases}
         \frac{(m-1)c_0\, p_c(X_0^*)^2}{\alpha-2} \, \frac{1}{M^{\alpha-2}}, & \hbox{\it if $2<\alpha\le 4$},\\
         \frac{1}{(m-1) c_0} \, \frac{1}{\log M}, & \hbox{\it if $\alpha=2$}.
     \end{cases}
     \label{equiv:p_M-p_c}
 \end{equation}

\end{lemma}

\medskip

\noindent {\it Proof of Lemma \ref{l:p_M-p_c}.} We first treat the case $2<\alpha\le 4$. Since $a_M=0$, \eqref{1/pc-1/pM} becomes $\frac{1}{p_c(X_0^*)} - \frac{1}{p_M} = \P(X_0^*>M) + \E \{ [(m-1) X_0^* -1] m^{X_0^*} \, {\bf 1}_{\{ X_0^* > M\} } \}$; using the asymptotics of $X_0^*$ given in \eqref{initial_law}, we obtain
\begin{equation}
    \frac{1}{p_c(X_0^*)} - \frac{1}{p_M}
    \, \sim \,
    \frac{(m-1)c_0}{\alpha-2} \, \frac{1}{M^{\alpha-2}} \, ,
    \qquad
    M\to \infty\, .
    \label{1/p_c-1/p_M}
\end{equation}

\noindent This yields \eqref{equiv:p_M-p_c} in case $2<\alpha \le 4$.

We now treat the case $\alpha=2$. By the asymptotics of $X_0^*$ given in \eqref{initial_law}, we have, in case $\alpha=2$,
$$
\E \{ [(m-1) X_0^* -1] m^{X_0^*} \, {\bf 1}_{\{ X_0^* \le M\} } \}
\, \sim \,
(m-1) c_0 \log M \, ,
\qquad
M\to \infty \, .
$$

\noindent In view of \eqref{1/p_M} (and using the fact $a_M=0$, as well as the trivial inequality $\P(X_0^*>M) \le 1$), this implies
\begin{equation}
    \frac{1}{p_M}
    \, \sim \,
    (m-1) c_0 \, \log M \, ,
    \qquad
    M\to \infty\, ,
    \label{1/p_M-ub}
\end{equation}

\noindent which, in turn, yields \eqref{equiv:p_M-p_c} in case $\alpha =2$.\qed

\bigskip

We have now all the ingredients for the proof of the upper bound in Theorems \ref{t:2<alpha<4} and \ref{t:alpha=2}.

Let $(Y_n^{(M)}, \, n\ge 0)$ be a recursive system with initial distribution $Y_0^{(M)}$. Then the system is critical; in fact, \eqref{p_M} is a rewriting of \eqref{p_c}. By Theorem A (see the introduction), the free energy of the system is $0$. It follows from the first inequality in \eqref{F_encadrement} that $\E_{p_M} ( Y_n^{(M)} ) \le \frac{1}{m-1}$ for all $n\ge 0$.

Consider also the system $(X_n, \, n\ge 0)$ under $\P_{p_M}$, i.e., with initial law $P_{X_0} = p_M \, P_{X_0^*} + (1-p_M) \delta_0$. Since $X_0 \ge Y_0^{(M)}$ $\P_{p_M}$-a.s., we are entitled to apply Lemma \ref{l:comparison} to see that for all $n\ge 0$,
\begin{equation}
    0 \le \E_{p_M}(X_n)-\E_{p_M}(Y_n^{(M)}) \le m^n \, [\E_{p_M}(X_0)-\E_{p_M}(Y_0^{(M)})] \, .
    \label{EX-EY:truncation}
\end{equation}

\noindent We take $n=M$. On the right-hand side,  note that $\E_{p_M}(X_0)-\E_{p_M}(Y_0^{(M)}) = p_M \, \E(X_0^* \, {\bf 1}_{\{ X_0^*>M\} }) \le \E(X_0^* \, {\bf 1}_{\{ X_0^*>M\} })$, so
$$
m^M \, [\E_{p_M}(X_0)-\E_{p_M}(Y_0^{(M)})]
\le
m^M \sum_{k=M+1}^\infty k \, \P(X_0^* =k) \, ,
$$

\noindent which tends to $0$ as $M\to \infty$ (by \eqref{initial_law}). Consequently, \eqref{EX-EY:truncation} yields that for all $M$ sufficiently large (say $M\ge M_0$), we have $0 \le \E_{p_M}(X_M)-\E_{p_M}(Y_M^{(M)}) \le 1$. We have already seen that $\E_{p_M}(Y_M^{(M)}) \le \frac{1}{m-1}$. Hence $\E_{p_M}(X_M) \le 1+\frac{1}{m-1} = \frac{m}{m-1}$ (for $M\ge M_0$). By the second inequality in \eqref{F_encadrement}, we obtain that the free energy of the system $(X_n, \, n\ge 0)$ with $p=p_M$ satisfies
\begin{equation}
    F_\infty (p_M)
    \le
    \frac{\E_{p_M}(X_M)}{m^M}
    \le
    \frac{m}{m-1} \, \ee^{-M \log m} \, ,
    \qquad
    \forall M\ge M_0\, .
    \label{F_infty(p_M)<}
\end{equation}

Let $p\in [p_{M+1}, \, p_M]$. When $p$ is sufficiently close to $p_c(X_0^*)$ (which is equivalent to saying that $M$ is sufficiently large), we have, by \eqref{F_infty(p_M)<},
\begin{equation}
    F_\infty(p)
    \le
    F_\infty(p_M)
    \le
    \frac{m}{m-1} \, \ee^{-M \log m} \, .
    \label{ub_finalstep}
\end{equation}

In case $2<\alpha\le 4$, we deduce from Lemma \ref{l:p_M-p_c}     that $\frac{m}{m-1} \, \ee^{-M \log m} = \exp ( - \frac{c_{47} + o(1)}{(p_M - p_c(X_0^*))^{1/(\alpha-2)}})$, where $c_{47}:= (\frac{(m-1)c_0\, p_c(X_0^*)^2}{\alpha-2})^{1/(\alpha-2)} \log m$. Note that as $p\downarrow p_c(X_0^*)$, $\frac{1}{(p_M - p_c(X_0^*))^{1/(\alpha-2)}} = \frac{1 + o(1)}{(p_{M+1} - p_c(X_0^*))^{1/(\alpha-2)}} \ge \frac{1 + o(1)}{(p - p_c(X_0^*))^{1/(\alpha-2)}}.$ Hence
\begin{equation}
    F_\infty(p)
    \le
    \exp \Big( - \frac{c_{47} + o(1)}{(p - p_c(X_0^*))^{1/(\alpha-2)}} \Big),
    \qquad
    p\downarrow p_c(X_0^*) \, .
    \label{free-energy_2<alpha<4:ub}
\end{equation}

\noindent This yields the upper bound in Theorem \ref{t:2<alpha<4}.\footnote{The proof of \eqref{free-energy_2<alpha<4:ub} does not require $\alpha\le 4$. In other words, \eqref{free-energy_2<alpha<4:ub} holds for all $\alpha> 2$; however, the exponent $\frac{1}{\alpha-2}$ is not optimal when $\alpha>4$, as we have seen in Theorem \ref{t:main}.}

In case $\alpha=2$, Lemma \ref{l:p_M-p_c} says that $\frac{m}{m-1} \, \ee^{-M \log m} = \exp ( - \ee^{(1+o(1))/[(m-1)c_0\, p]})$, $p\downarrow 0$. The upper bound in Theorem \ref{t:alpha=2} follows immediately from \eqref{ub_finalstep}.\qed

\subsection{Proof of Theorems \ref{t:2<alpha<4} and \ref{t:alpha=2}: lower bound}
\label{subs:light_tail_lb}

To prove the lower bound, here is our choice of $\varrho$ in \eqref{X0M}: for $\alpha=2$, $\varrho$ is simply $0$; for $2<\alpha \le 4$, $\varrho$ is chosen to be independent of $X_0^*$, with $\P(\varrho = 0) = \frac{1}{M^{\alpha-2}} = 1- \P(\varrho =\ell_0)$.

\medskip

\begin{lemma}
\label{l:p_M-p_c-lb}

 Assume $\P(X_0^* = k) \sim c_0 \, m^{-k} k^{-\alpha}$, $k\to \infty$, for some $0<c_0<\infty$ and $2\le \alpha \le 4$. When $M\to \infty$,
 \begin{equation}
     p_M-p_c(X_0^*)
     \, \sim \,
     \begin{cases}
         \frac{c_{48}}{M^{\alpha-2}}, & \hbox{\it if $2<\alpha\le 4$},\\
         \frac{1}{(m-1) c_0} \, \frac{1}{\log M}, & \hbox{\it if $\alpha=2$},
     \end{cases}
     \label{equiv:p_M-p_c-lb}
 \end{equation}
 where, in case $2<\alpha\le 4$, $c_{48} := [\, p_c(X_0^*)]^2 \, c_{49}>0$, with
 \begin{equation}
     c_{49}
     :=
     \frac{(m-1)c_0}{\alpha-2}
     +
     \{ [(m-1)\ell_0 -1]m^{\ell_0} +1 \} \, \P(X_0^* = \ell_0) \, .
     \label{cst}
 \end{equation}

\end{lemma}

\medskip

\noindent {\it Proof of Lemma \ref{l:p_M-p_c-lb}.} The proof is as in the proof of Lemma \ref{l:p_M-p_c}, except that this time, we need to evaluate $a_M$ with our new choice of $\varrho$. Let us only point out the difference.

In case $2<\alpha\le 4$, we have, by definition of $a_M$ in \eqref{a_M},
$$
a_M
=
\E \{ [(m-1)\ell_0 -1]m^{\ell_0} - [(m-1)\varrho -1]m^\varrho \} \, \P(X_0^* = \ell_0)
=
\frac{c_{50}}{M^{\alpha-2}}\, ,
$$

\noindent where $c_{50} := \{ [(m-1)\ell_0 -1]m^{\ell_0} +1 \} \, \P(X_0^* = \ell_0) >0$. Then \eqref{1/p_c-1/p_M} in the proof of Lemma \ref{l:p_M-p_c} becomes
$$
\frac{1}{p_c(X_0^*)} - \frac{1}{p_M}
\, \sim \,
\frac{(m-1)c_0}{\alpha-2} \, \frac{1}{M^{\alpha-2}}
+
a_M
=
\frac{c_{49}}{M^{\alpha-2}}\, ,
\qquad
M\to \infty\, ,
$$

\noindent with $c_{49}>0$ being in \eqref{cst}. This yields Lemma \ref{l:p_M-p_c-lb} in case $2<\alpha\le 4$.

Assume now $\alpha=2$ (so $p_c(X_0^*)=0$). We have chosen $\varrho=0$ in this case. By definition of $a_M$ in \eqref{a_M},
$$
a_M
=
\{ [(m-1)\ell_0 -1]m^{\ell_0} +1 \} \, \P(X_0^* = \ell_0)
=
c_{50} > 0 \, .
$$

\noindent Then \eqref{1/p_M-ub} in the proof of Lemma \ref{l:p_M-p_c} becomes
$$
\frac{1}{p_M}
\, \sim \,
(m-1) c_0 \, (\log M)
-
c_{50}
\, \sim \,
(m-1) c_0 \, \log M \, ,
\qquad
M\to \infty\, .
$$

\noindent This yields Lemma \ref{l:p_M-p_c-lb} in case $\alpha=2$.\qed

\bigskip

Recall from \eqref{alpha} that a system with initial value $X_0$ is said to be $\beta$-regular with coefficient $\chi$ if for all integers $j\ge 1$,
$$
\Xi_j(X_0)
:=
\E \Big[ (X_0 \wedge j)^2 \, [(m-1)X_0-1] \, m^{X_0} \Big]
\ge
\chi\, \min\{\E(X_0^3 m^{X_0}), \, j^\beta\} \, .
$$

\medskip

\begin{lemma}
 \label{e:chi_lower_bound}

 Let $\alpha\in [2, \, 4]$.
 There exists a constant $c_{51}\in (0, \, 1]$ whose value depends only on $m$, $\alpha$, and the distribution of $X_0^*$, such that for all sufficiently large $M$, the critical system $(Y_n^{(M)}, \, n\ge 0)$ is $(4-\alpha)$-regular with coefficient $c_{51}$.

\end{lemma}

\medskip

\noindent {\it Proof.} When $\alpha=4$, Lemma \ref{A_4P(X_0=0)} says that $(Y_n^{(M)}, \, n\ge 0)$ is $0$-regular with coefficient $1-p_M$, which goes to $1-p_c(X_0^*)>0$ when $M\to \infty$, so the result follows immediately.

Let $\alpha\in [2, \, 4)$ now. By assumption \eqref{initial_law}, $\P (X_0^* = k) \, \sim \, c_0 \, m^{-k} k^{-\alpha}$, $k\to \infty$, so there exist constants $c_{52}>0$ and $c_{53}>0$, and an integer $k_0\ge 2$ (for further use, we take $k_0$ to satisfy also $k_0\ge \frac{2}{m-1}$), all depending only on $m$, $\alpha$, and the law of $X_0^*$, such that for $M\ge k_0$,
\begin{eqnarray}
    \P (Y_0^{(M)} = k)
 &\le& c_{52} \, m^{-k} k^{-\alpha}, \qquad \forall 1\le k\le M,
    \label{c820}
    \\
    \P (Y_0^{(M)} = k)
 &\ge& c_{53} \, m^{-k} k^{-\alpha}, \qquad \forall k_0\le k\le M.
    \label{c802}
\end{eqnarray}

\noindent So for all integer $M\ge 2k_0$ and real number $b\in [k_0, \, \frac{M}{2}]$, we have, by \eqref{c802},
\begin{eqnarray}
    \E [ m^{Y_0^{(M)}} Y_0^{(M)} \, {\bf 1}_{\{ b\le Y_0^{(M)}\le 2b\} } ]
 &\ge& \sum_{k\in [b, \, 2b] \cap \z} m^k k \, c_{53} m^{-k} k^{-\alpha}
   \nonumber
    \\
    &\ge&
    c_{53}\sum_{k\in [b, \, 2b] \cap \z} (2b)^{1-\alpha}
    \nonumber
    \\
 &\ge& c_{54}\, b^{2-\alpha} \, ,
    \label{truncated_moment_lb}
\end{eqnarray}

\noindent with $c_{54} := 2^{-\alpha} \, c_{53}$. Furthermore, recalling that we only care about the case $\alpha\in [2, \, 4)$, we have, by \eqref{c820},
\begin{equation}
    \E((Y_0^{(M)})^3 m^{Y_0^{(M)}})
    \le
    \sum_{k=1}^M k^3 m^{k} \, c_{52} \, m^{-k} k^{-\alpha}
    \le
    c_{55} \, M^{4-\alpha} \, ,
    \label{Lambda_truncated}
\end{equation}

\noindent where $c_{55} := \frac{2^{4-\alpha}\, c_{52}}{4-\alpha}$ (because $\sum_{k=1}^M k^{3-\alpha} \le \int_0^{M+1} x^{3-\alpha} \d x = \frac{(M+1)^{4-\alpha}}{4-\alpha} = \frac{M^{4-\alpha}}{4-\alpha} \, (1+\frac1M)^{4-\alpha} \le \frac{M^{4-\alpha}}{4-\alpha} \, 2^{4-\alpha}$).

Let $j\ge 1$ be an integer. We distinguish two possible situations.

{\it First situation: $M\ge 2k_0j$}. We use
\begin{eqnarray*}
    \Xi_j(Y_0^{(M)})
 &:=& \E [ m^{Y_0^{(M)}} [(m-1)Y_0^{(M)}-1] (Y_0^{(M)} \wedge j)^2]
    \\
 &\ge& \E [ m^{Y_0^{(M)}} [(m-1)Y_0^{(M)}-1] (Y_0^{(M)} \wedge j)^2 \, {\bf 1}_{\{ k_0j\le Y_0^{(M)}\le 2k_0j\} } ] \, .
\end{eqnarray*}

\noindent On the event $\{ k_0j\le Y_0^{(M)}\le 2k_0j\}$ (so $Y_0^{(M)} \ge \frac{2}{m-1}$), we have $(m-1)Y_0^{(M)}-1 \ge \frac{m-1}{2} \, Y_0^{(M)}$ and $Y_0^{(M)} \wedge j = j$, so
$$
\Xi_j(Y_0^{(M)})
\ge
\frac{m-1}{2} \, j^2 \, \E [ m^{Y_0^{(M)}} Y_0^{(M)} \, {\bf 1}_{\{ k_0j\le Y_0^{(M)}\le 2k_0j\} } ]
\ge
\frac{m-1}{2}\, c_{54}\, k_0^{2-\alpha}\, j^{4-\alpha} \, ,
$$

\noindent the last inequality being a consequence of \eqref{truncated_moment_lb}.

{\it Second (and last) situation: $M<2k_0j$} (and $M\ge 2k_0$). This time, we start with
$$
\Xi_j(Y_0^{(M)})
\ge
\E [ m^{Y_0^{(M)}} [(m-1)Y_0^{(M)}-1] (Y_0^{(M)} \wedge j)^2 \, {\bf 1}_{\{ \frac{M}{2} \le Y_0^{(M)}\le M\} } ] \, .
$$

\noindent On the event $\{ \frac{M}{2} \le Y_0^{(M)}\le M\}$ (so $Y_0^{(M)} \ge \frac{2}{m-1}$), we have again $(m-1)Y_0^{(M)}-1 \ge \frac{m-1}{2} \, Y_0^{(M)}$, whereas $(Y_0^{(M)} \wedge j) \ge \frac{M}{2k_0}$. Consequently,
\begin{eqnarray*}
\Xi_j(Y_0^{(M)})
&\ge&
\frac{m-1}{2} \, \frac{M^2}{4k_0^2} \, \E [ m^{Y_0^{(M)}} Y_0^{(M)} \, {\bf 1}_{\{ \frac{M}{2} \le Y_0^{(M)}\le M\} } ]
\\
&\ge&
\frac{m-1}{2} \, \frac{M^2}{4k_0^2} \, c_{54}\, (\frac{M}{2})^{2-\alpha} \, ,
\end{eqnarray*}

\noindent the last inequality following again from \eqref{truncated_moment_lb}. By \eqref{Lambda_truncated}, $\E((Y_0^{(M)})^3 m^{Y_0^{(M)}}) \le c_{55} \, M^{4-\alpha}$, so $\Xi_j(Y_0^{(M)}) \ge c_{56} \, \E((Y_0^{(M)})^3 m^{Y_0^{(M)}})$, with $c_{56} := \frac{m-1}{2} \, \frac{c_{54}}{2^{4-\alpha} k_0^2c_{55}}$.

In both situations, we see that for all $M\ge 2k_0$ and all integers $j\ge 1$,
$$
\Xi_j(Y_0^{(M)})
\ge
c_{57} \, \min\{ \E((Y_0^{(M)})^3 m^{Y_0^{(M)}}), \, j^{4-\alpha}\} \, ,
$$

\noindent with $c_{57} := \min\{ \frac{m-1}{2} c_{54}k_0^{2-\alpha}, \, c_{56}\}$. In other words,  for all $M\ge 2k_0$, the system $(Y_n^{(M)}, \, n\ge 0)$ is $(4-\alpha)$-regular with coefficient $\min\{ c_{57}, \, 1\}$. Lemma \ref{e:chi_lower_bound} is proved.\qed

\bigskip

We now proceed to the proof of the lower bound in Theorems \ref{t:2<alpha<4} and \ref{t:alpha=2}. Again, we consider the system $(X_n, \, n\ge 0)$ under $\P_{\! p_M}$, i.e., with initial law $P_{X_0} = p_M \, P_{X_0^*} + (1-p_M) \delta_0$, and couple it with the critical system $(Y_n^{(M)}, \, n\ge 0)$.

By definition, $\P_{\! p_M} (X_0=k) \ge \P(Y_0^{(M)}=k)$ for all integers $k\ge 1$. We intend to apply Theorem \ref{t:2nd_crossing} to $(X_n, \, n\ge 0)$ under $\P_{\! p_M}$, so let us estimate $\E [(Y_0^{(M)})^3 m^{Y_0^{(M)}}]$ and $\E_{p_M} (X_0) - \E (Y_0^{(M)})$.

Since $X_0^{(M)} \le X_0^*\, {\bf 1}_{\{ X_0^* \le M\} }$, we have
$$
\E [(Y_0^{(M)})^3 m^{Y_0^{(M)}}]
=
p_M \, \E [(X_0^{(M)})^3 m^{X_0^{(M)}}]
\le
p_M \, \E [(X_0^*)^3 m^{X_0^*}\, {\bf 1}_{\{ X_0^* \le M\} }] \, .
$$

\noindent By the tail behaviour of $X_0^*$ given in \eqref{initial_law}, we know that for $M\to \infty$, $\E [(X_0^*)^3 m^{X_0^*}\, {\bf 1}_{\{ X_0^* \le M\} }]$ is equivalent to $\frac{c_0}{4-\alpha} \, M^{4-\alpha}$ if $2\le \alpha <4$, and to $c_0 \log M$ if $\alpha=4$. So in both cases,
\begin{equation}
    \E [(Y_0^{(M)})^3 m^{Y_0^{(M)}}]
    \le
    p_M \, M^{4-\alpha+o(1)} \, ,
    \qquad
    M\to \infty\, .
    \label{truncating_argument_1}
\end{equation}

Let us now estimate $\E_{p_M} (X_0) - \E (Y_0^{(M)})$. By definition,
\begin{eqnarray*}
    \E_{p_M} (X_0) - \E (Y_0^{(M)})
 &=& p_M \, [\E (X_0^*) - \E (X_0^{(M)})]
    \\
 &\ge& p_M \, [\E (X_0^* \, {\bf 1}_{\{ X_0^* \le M\} }) - \E (X_0^{(M)})]
    \\
 &=& p_M \, \E [(\ell_0 - \varrho)\, {\bf 1}_{\{ X_0^*=\ell_0\} }]
    \\
 &=& p_M \, [\ell_0 -\E(\varrho)] \, \P (X_0^*=\ell_0)\, .
\end{eqnarray*}

\noindent For $2<\alpha\le 4$, $\P(\varrho = 0) = \frac{1}{M^{\alpha-2}} = 1- \P(\varrho =\ell_0)$, so $\ell_0 -\E(\varrho) = \frac{\ell_0}{M^{\alpha-2}}$. For $\alpha=2$, $\varrho=0$ so $\ell_0 -\E(\varrho) = \ell_0$. In both cases, $\ell_0 -\E(\varrho) = \frac{\ell_0}{M^{\alpha-2}}$, hence
\begin{equation}
    \E_{p_M} (X_0) - \E (Y_0^{(M)})
    \ge
    \frac{\ell_0\, \P (X_0^*=\ell_0)}{M^{\alpha-2}} \, p_M \, .
    \label{truncating_argument_2}
\end{equation}

Lemma \ref{e:chi_lower_bound} says the system $(Y_n^{(M)}, \, n\ge 0)$ is $(4-\alpha)$-regular in the sense of \eqref{alpha}, with coefficient $\chi := c_{51}$. Now that \eqref{truncating_argument_1} and \eqref{truncating_argument_2} are established, and since $\P(Y_0^{(M)} =0) \ge 1-p_M \ge \frac{1-p_c(X_0^*)}{2}>0$ (for all sufficiently large $M$), we are entitled to apply Theorem \ref{t:2nd_crossing} to $(X_n, \, n\ge 0)$ under $\P_{\! p_M}$, to see that for $M\to \infty$, the free energy of $(X_n, \, n\ge 0)$ under $\P_{\! p_M}$ satisfies
$$
F_\infty (p_M)
\ge
\exp \Big( - M^{1+o(1)} [\log (M^{2+o(1)})]^{\max\{(\alpha/2), \, 2\}} \Big)
=
\exp \Big( - M^{1+o(1)} \Big)\, .
$$

\noindent By Lemma \ref{l:p_M-p_c-lb}, for $M\to \infty$, $M^{1+o(1)} = \frac{1}{(p_M-p_c(X_0^*))^{1/(\alpha-2)+o(1)}}$ if $2<\alpha\le 4$, and $M^{1+o(1)} = \exp( \frac{1+o(1)}{(m-1)c_0} \, \frac{1}{p_M})$ if $\alpha=2$. This yields the lower bound in Theorems \ref{t:2<alpha<4} and \ref{t:alpha=2} along the sequence $p=p_M$. The monotonicity of $p\mapsto F_\infty(p)$ and the properties $\frac{p_M-p_c(X_0^*)}{p_{M+1}-p_c(X_0^*)} \to 1$ (if $2<\alpha\le 4$) and $\frac{p_M}{p_{M+1}} \to 1$ (if $\alpha=2$; both properties are consequences of Lemma \ref{l:p_M-p_c-lb}) allow us, exactly as in the proof of the upper bound in Theorems \ref{t:2<alpha<4} and \ref{t:alpha=2}, to get rid of the restriction to the sequence $(p_M)$.\qed

\bigskip

\noindent {\bf Acknowledgements.} We are grateful to two anonymous referees, whose insightful comments have led to improvements in the presentation of the paper.

\end{document}